\documentclass[12pt]{amsart} 

\usepackage{amsmath,amssymb,color,amsthm,epsfig,extarrows,tikz,enumitem,
hyperref,url}

\theoremstyle{definition}
\newtheorem{prop}{Proposition}[section]
\newtheorem{algor}[prop]{Algorithm} 
\newtheorem{thm}[prop]{Theorem}
\newtheorem{rem}[prop]{Remark}
\newtheorem{lem}[prop]{Lemma}
\newtheorem{dfn}[prop]{Definition}

\newtheorem{mahler}{Mahler's Theorem}

\newtheorem{hensel}{Hensel's Lemma}

\newtheorem{ex}[prop]{Example}
\newtheorem{cor}[prop]{Corollary}

\newtheorem{not*}[prop]{Notation} 

\newcommand{\F}{\mathbb{F}}
\newcommand{\N}{\mathbb{N}}
\newcommand{\Q}{\mathbb{Q}}
\newcommand{\R}{\mathbb{R}}
\newcommand{\C}{\mathbb{C}}

\newcommand{\Z}{\mathbb{Z}}

\newcommand{\cT}{\mathcal{T}}

\newcommand{\np}{\mathbf{NP}}
\newcommand{\zpp}{\mathbf{ZPP}}

\newcommand{\norm}[1]{\left|#1\right|_p}
\newcommand{\ord}{\operatorname{ord}}
\newcommand{\newt}{\operatorname{Newt}}
\newcommand{\anewt}{\operatorname{Newt}_\infty}

\newcommand{\tf}{\tilde{f}}

\renewcommand{\qed}{$\blacksquare$}
\newcommand{\dia}{$\diamond$}
\newcommand{\eps}{\varepsilon}

\newcommand{\floor}[1]{\left\lfloor #1 \right\rfloor}
\newcommand{\ceil}[1]{\left\lceil #1 \right\rceil}

\title[A Complexity Chasm for Sparse Solving Over $\Q_p$]{\mbox{}\\ 
\vspace{-1.3in}A Complexity Chasm for Solving Univariate Sparse Polynomial 
Equations Over $p$-adic Fields}

\author{J.\ Maurice Rojas}
\thanks{Partially supported by NSF grants CCF-1900881 and 
CCF-1409020. This paper contains an Appendix missing from the 
proceedings version \cite{issacversion}, as well as some   
corrections and improvements.}
\email{jmauricerojas@gmail.com}
\author{Yuyu Zhu}
\email{yuyu.zhu1213@gmail.com}
\address{Texas A\&{}M University, TAMU 3368, College Station, Texas \ 
77843-3368 }

\begin{document} 

\begin{abstract}
We reveal a complexity chasm, separating the trinomial and 
tetranomial cases, for solving univariate sparse polynomial equations 
over certain local fields. First, for any fixed field 
$K\!\in\!\{\Q_2,\Q_3,\Q_5,\ldots\}$, we prove that any polynomial 
$f\!\in\!\Z[x]$ with exactly $3$ monomial terms, 
degree $d$, and all coefficients having absolute value at most $H$, 
can be solved over $K$ within deterministic time $\log^{O(1)}(dH)$ in the 
classical Turing model. (The best previous algorithms were of complexity 
exponential in $\log d$, even for just counting roots in $\Q_p$.) In 
particular, our algorithm generates approximations 
in $\Q$ with bit-length $\log^{O(1)}(dH)$ to all the roots of $f$ in $K$, 
and these approximations converge quadratically under Newton iteration. 
On the other hand, we give a unified family of {\em tetra}nomials requiring  
$\Omega(d\log H)$ digits to distinguish the base-$p$ 
expansions of their roots in $K$. 
\end{abstract}


\maketitle

\vspace{-.4in} 
\section{Introduction} 
Solving polynomial equations over the $p$-adic rational numbers $\Q_p$ 
underlies many classical questions in number theory, and is close to numerous 
applications in cryptography, coding theory, and computational number theory. 
Furthermore, the complexity of 
solving {\em structured} equations --- such as those with a fixed number of 
monomial terms or invariance with respect to a group action --- arises 
naturally in many computational geometric applications and 
is closely related to a deeper understanding of circuit complexity 
(see, e.g., \cite{koiranrealtau}). 
So we will classify when it is possible to separate and approximate 
roots in $\Q_p$ in deterministic polynomial-time. 

Recall that thanks to 17th century work of Descartes, and 20th century 
work of Lenstra \cite{len99} and Poonen \cite{poonen},    
it is known that univariate polynomials with exactly $t$ monomial 
terms have at most $t^{O(1)}$ roots in a fixed field $K$ {\em only} 
when $K$ is $\R$ or a finite algebraic extension of $\Q_p$ for some prime 
$p\!\in\!\N$. 
We'll use $|\cdot|_p$ (resp.\ $|\cdot|$) for the absolute value on the 
$p$-adic complex numbers $\C_p$ \cite{Rob00} 
normalized so that $|p|_p\!=\!\frac{1}{p}$ (resp.\ the standard absolute 
value on $\C$). Recall also that for any function $f$ analytic on $K$, the 
corresponding {\em Newton endomorphism} is $N_f(z):=z-\frac{f(z)}{f'(z)}$, and 
the corresponding sequence of {\em Newton iterates} of a {\em start-point} 
$z_0\!\in\!K$ is the sequence $(z_i)^\infty_{i=0}$ where 
$z_{i+1}\!:=\!N_f(z_i)$ for all $i\!\geq\!0$. Finally, we call any 
polynomial in $\Z[x_1,\ldots,x_n]$ having exactly $t$ terms in its monomial 
term expansion an {\em $n$-variate $t$-nomial}. We will often use 
$x$ in place of $x_1$. 

Our first main result is that we can efficiently count the roots of 
trinomials in $\Q_p$, {\em and} find succinct start-points in $\Q$ 
under which Newton iteration converges quickly to all the roots in $\Q_p$. 
We use $\#S$ for the cardinality of a set $S$.  

\begin{thm} 
\label{thm:big} {\em 
Suppose $K\!=\!\Q_p$ for some {\em fixed}\footnote{We clarify 
the dependence of our complexity bounds on $p$ in Section 
\ref{sec:trinosolqp}.} prime $p\!\in\!\N$. Then 
for any input trinomial $f\!\in\!\Z[x]$ with 
degree $d$ and all coefficients of (Archimedean) absolute value
$\leq\!H$, we can find in deterministic time 
$O\!\left(\log^{16}(dH)\log\log(dH)\right)$ a   
set $\{\frac{a_1}{b_1},\ldots,\frac{a_m}{b_m}
\}\subset\!\Q$ of cardinality $m\!=\!m(K,f)$ such that:\\  
\mbox{}\hspace{.5cm}1. For all $j$ we have $a_j\!\neq\!0 
\Longrightarrow \log|a_j|,\log|b_j|=O(\log^8(dH))$. \\ 
\mbox{}\hspace{.5cm}2. There is a $\mu\!=\!\mu(d,H)\!>\!1$ 
such that $z_0\!:=\!a_j/b_j$ implies that $f$ 
has a root $\zeta_j\in K$ with\\   
\mbox{}\hspace{1cm}sequence of Newton iterates satisfying 
$|z_i-\zeta_j|_p\!\leq\!\mu^{-2^{i-1}}|z_0-\zeta_j|_p$ 
for all $i,j\!\geq\!1$.\\  
\mbox{}\hspace{.5cm}3. $m\!=\!\#\{\zeta_1,\ldots,\zeta_m\}$ is exactly the 
number of roots of $f$ in $K$.} 
\end{thm} 

\noindent 
We prove Theorem \ref{thm:big} in Section \ref{sec:trinosolqp}, 
via Algorithm \ref{algor:trinosolqp} there. 
(An analogue of Theorem \ref{thm:big} in fact holds for $K\!=\!\R$ as well, 
and will be presented in a sequel to  
this paper.)  We will call the convergence condition on $z_0$ above 
{\em being an approximate root (in the sense of Smale}\footnote{This 
terminology has only been applied over $\C$ with $\mu\!=\!2$ so far 
\cite{smale}, so we take 
the opportunity here to extend it to the $p$-adic rationals.}{\em \!)}, 
{\em with associated true root $\zeta_j$}. This type 
of convergence provides an efficient encoding of an approximation that 
can be quickly tuned to any desired accuracy. 

\begin{rem} {\em Defining the {\em input size} of a univariate polynomial 
$f(x)\!:=\!\sum^t_{i=1} c_i x^{a_i}\!\in\!\Z[x]$ as 
$\sum^t_{i=1}\log((|c_i|+2)(|a_i|+2))$ we see that Theorem \ref{thm:big} 
implies that one can solve univariate\linebreak 
\scalebox{.91}[1]{trinomial equations, over any {\em fixed} 
$p$-adic field, in deterministic time polynomial in the input size. \dia}} 
\end{rem} 
\begin{rem} {\em 
Efficiently solving univariate $t$-nomial equations over $K$ 
in the sense of Theorem \ref{thm:big} is easier for $t\!\leq\!2$: 
The case $t\!=\!1$ is clearly trivial (with $0$ the only possible root) 
while the case $(K,t)\!=\!(\R,2)$ is implicit in work 
on computer arithmetic from the 1970s (see, e.g., \cite{borwein}). 
We review the case $(K,t)\!=\!(\Q_p,2)$ with $p$ prime in Corollary 
\ref{cor:binomod} and Theorem \ref{thm:binoqp} of Section 2 below. \dia}  
\end{rem} 

Despite much work on factoring univariate polynomials over $\Q_p$ 
(see, e.g., \cite{CG00,GNP12,bns13,blq13}), all known general algorithms for 
solving (or even just counting the solutions of) arbitrary degree $d$ 
polynomial equations over $\Q_p$ have complexity exponential in $\log d$.  
So Theorem \ref{thm:big} presents a significant new speed-up, and improves an 
earlier complexity bound (membership in $\np$, for detecting roots in $\Q_p$) 
from \cite{airr}. We'll see in Section \ref{sec:trisepqp} how our speed-up 
depends on $p$-adic Diophantine approximation \cite{yu94}. 
Another key new ingredient in proving Theorem\linebreak 
\scalebox{.95}[1]{\ref{thm:big} is an efficient 
encoding of roots in $\Z/(p^k)$ from \cite{DMS19,krrz19} (with 
an important precursor in \cite{blq13}).}   

\subsection{Why is the Field Fixed?} 
Much as real algebraic geometry fixes the underlying field to 
$K\!=\!\R$ once and for all, our results focus on  
$K\!=\!\Q_p$ with $p$ fixed once and for all. In particular, while there are 
certainly number-theoretic
algorithms with deterministic complexity having 
dependence $(\log p)^{O(1)}$ on an input prime $p$, solving 
sparse polynomial equations in one variable 
over $\Q_p$ appears to have much larger complexity as a function of $p$. 
There are naive reasons, and subtle reasons, for this: 

\noindent 
{\bf R1.} {\em Whereas a binomial has at most $3$ roots in $\R$ (e.g., 
$x^3-x$), a binomial can have as many as $p$ roots in $\Q_p$ (e.g., 
$x^p-x$). Furthermore, trinomials have at most $5$, $7$, $9$, or 
$3p-2$\linebreak 
\scalebox{.95}[1]{roots in $K$, according as $K$ is $\R$, $\Q_2$ \cite{len99}, 
$\Q_3$ \cite{zhuthesis}, or $\Q_p$ with $p\!\geq\!5$ \cite{ak11}, and each 
bound is sharp.}} 

\noindent 
{\bf R2.} {\em Approximating square-roots of $p$-adic integers not divisible 
by $p$, within accuracy $1$, is equivalent to finding square-roots in 
the finite field $\F_p$. The latter problem is {\em still} not known to 
be doable in deterministic time polynomial in $\log p$, even though 
the decision version is doable in deterministic polynomial-time (see, e.g., 
\cite{bs,poonenzeta}).}  

In particular, even if one only wants to approximate just one or two  
roots in $\Q_p$, the minimal {\em currently} provable accuracy needed to 
decide if two approximations converge to the same root appears to have 
{\em quasi-linear} dependence on $p$. Interestingly, the truth of strong forms 
of the abc-Conjecture would imply a much smaller and practical dependence on 
$p$: See \cite{bakerabc} and Section \ref{sec:trisepqp} below. 

\subsection{The Separation Chasm at Four Terms} 
\label{sub:tetrasep} 
The $p$-adic rational roots of sparse polynomials can range from well-separated 
to (possibly) tightly spaced, already with just $4$ terms.  
\begin{thm} 
{\em \label{thm:tetra}
Consider the family of tetranomials \\ 
\mbox{}\hfill $\displaystyle{f_{d,\eps}(x):=x^d - \eps^{-2h}x^2 
+ 2\eps^{-(h+1)}x - \eps^{-2}}$\hfill\mbox{}\\ 
with $h\!\in\!\N$, $h\!\geq\!3$, and $d\!\in\!\{4,\ldots,\lfloor e^h 
\rfloor\}$ even. Let $H\!:=\!\max\{\eps^{\pm 2h}\}$.
Then $f_{d,\eps}$ has distinct roots $\zeta_1,\zeta_2\!\in\!K$ with 
$|\log|\zeta_1-\zeta_2|_p|$ or
$|\log|\zeta_1-\zeta_2||$ of order $\Omega(d\log H)$, according as 
$(K,\eps)\!=\!(\Q_p,p)$ or $(K,\eps)\!=\!(\R,1/2)$. In 
particular, the coefficients of $p^{2h}f_{d,p}$ (resp.\ $f_{d,\frac{1}{2}}$)  
all lie in $\Z$ and have $O(\log H)$ base-$p$ digits (resp.\ bits).}  
\end{thm}

\vspace{-.1cm} 
\noindent
We prove Theorem \ref{thm:tetra} in Section \ref{sec:tetra}. The 
special case $K\!=\!\R$ was derived earlier (in different notation) by 
Mignotte \cite{mig95}. (See also \cite{sag14}.) The case $K\!=\!\Q_p$ with $p$ prime appears to be new, and our proof unifies the Archimedean and 
non-Archimedean cases via tropical geometry.  
Note that Theorem \ref{thm:tetra} implies that the roots in $K$ of a 
tetranomial can be so close that 
one needs $\Omega(d\log H)$ many digits to distinguish their base-$p$ 
expansions in the worst case.    

Mignotte used the tetranomial $f_{d,1/2}$ in \cite{mig95} 
to show that an earlier root separation bound of Mahler 
\cite{mah64}, for {\em arbitrary} degree $d$ polynomials in $\Z[x]$, 
is asymptotically near-optimal. We recall the following paraphrased version: 

\begin{mahler} {\em 
Suppose $f\!\in\!\Z[x]$ has degree $d$, all coefficients of (Archimedean) 
absolute value at most $H$, and is irreducible in $\Z[x]$. Let 
$\zeta_1,\zeta_2\!\in\!\C$ be distinct roots of $f$. Then 
$|\log|\zeta_1-\zeta_2||\!=\!O(d\log(dH))$. \qed} 
\end{mahler} 

Our new algorithmic results are enabled by our third and final main
result: Mahler's bound can be dramatically 
improved for trinomials.   
\begin{thm} {\em 
\label{thm:tri} 
Suppose $p$ is prime and $f\!\in\!\Z[x]$ is square-free, 
has exactly $3$ monomial terms, degree $d$, and all coefficients of 
(Archimedean) absolute value at most $H$. Let $\zeta_1,\zeta_2\!\in\!\C_p$ be 
distinct roots of $f$. Then $|\log|\zeta_1-\zeta_2|_p|\!=\!
O\!\left(\frac{p}{\log^2 p}\log(d)\log^2(dH+p)\log\log(dH+p)\right)$.}    
\end{thm} 

\vspace{-.1cm} 
\noindent 
We prove Theorem \ref{thm:tri} in Section \ref{sec:trisepqp}. 
Theorem \ref{thm:tri} is in fact a $p$-adic analogue of a separation 
bound of Koiran for roots in $\R$ \cite{koiransep}. Even sharper bounds 
can be derived for binomials: We review these
bounds in Section \ref{sec:bisep}.

\subsection{Previous Complexity and Sparsity Results} 
Deciding the existence of roots 
over $\Q_p$ for univariate polynomials with an {\em arbitrary} number 
of monomial terms is already $\np$-hard with respect 
to randomized ($\zpp$, a.k.a.
Las Vegas) reductions \cite{airr}. 
On the other hand, detecting roots over $\Q_p$ for $n$-variate $(n+1)$-nomials 
is known to be doable in $\np$ \cite{airr}. Speeding this up to 
polynomial-time, even for $n\!=\!2$ and fixed $p$, hinges upon detecting roots 
in $(\Z/(p^k))^2$ for bivariate trinomials of degree $d$ 
in time $(k+\log d)^{O(1)}$. The latter problem remains open, but some 
progress has been made in author Zhu's Ph.D.\ thesis \cite{zhuthesis}. 

On a related note, counting points on trinomial curves over the prime 
fields $\F_p$ in time $(\log(pd))^{O(1)}$ remains a challenging open question. 
Useful quantitative estimates in this direction were derived in  
\cite{huavandiver} and revisited via real quadratic optimization in 
\cite{avendanomorales}.  

\section{Background} \label{sec:back} 
\subsection{Newton Polygons and Newton Iteration: Archimedean and 
Non-Archimedean}  
Definitive sources for $p$-adic arithmetic and analysis include 
\cite{serre,schikhof,Rob00}. 
We use $\ord_p : \C_p \longrightarrow \Q$ for the standard $p$-adic 
valuation on $\C_p$, normalized so that $\ord_p p\!=\!1$. 
The {\em most significant ($p$-adic)} {\em digit} of 
$\sum^\infty_{j=s}a_j p^j\!\in\!\Q_p$ is $a_s$, assuming the 
$a_j\!\in\!\{0,\ldots,p-1\}$ and $a_s\!\neq\!0$.  

The notion of Newton polygon goes back to 17th century work of Newton 
on Puiseux series solutions to polynomial equations 
\cite[pp.\ 126--127]{newtonbook}. 
We will need variants 
of this notion over $\C_p$ and $\C$. (See, e.g., \cite{wei63} for the 
$p$-adic case and \cite{ostrowskiarch,aknr} for the complex case.) 
\begin{dfn} {\em \label{dfn:newt} Suppose 
$f(x)\!:=\!\sum_{i=1}^{t}c_ix^{a_i}\!\in\!\Z[x]$ with 
$c_i\!\neq 0$ for all $i$ and $a_1\!<\!\cdots\!<\!a_t$. We then define 
the {\em $p$-adic Newton polygon}, $\newt_p(f)$ (resp.\ {\em 
Archimedean Newton polygon,} $\anewt(f)$) to be the convex hull of the set of 
points $\{(a_i,\ord_pc_i)\; | \;  i\!\in\!\{1,\ldots,t\}\}$ (resp.\ the convex 
hull of $\{(a_i,-\log|c_i|)\; | \;  i\!\in\!\{1,\ldots,t\}\}$). We call an edge 
$E$ of a polygon in $\R^2$ {\em lower} if and only if $E$ has an inner 
normal with positive last coordinate. We also define the 
{\em horizontal length} of a line segment $E$ connecting $(r,s)$ and $(u,v)$ 
to be $\lambda(E)\!:=\!|u-r|$. \dia } 
\end{dfn}
\begin{ex} {\em \label{ex:newts} 
Following the notation of Theorem \ref{thm:tetra}, we set 
$h\!=\!3$ and illustrate $\newt_p\left(f_{5,p}\right)$ (for $p$ odd) and 
$\newt_\infty(f_{5,1/2})$ below: \\  
\mbox{}\hfill\scalebox{1}[.7]{\epsfig{file=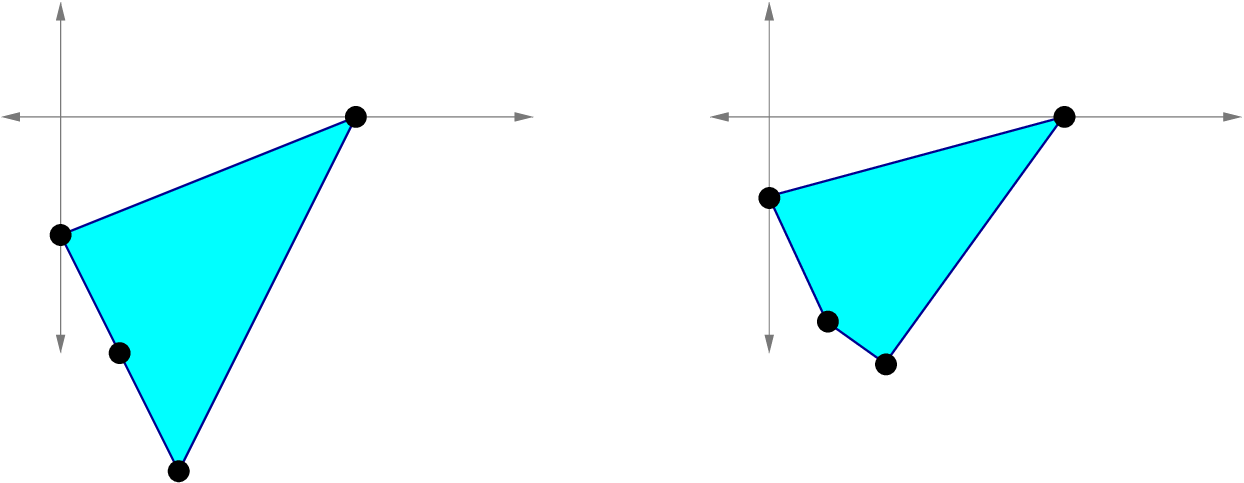,height=1.5in,
clip=}}\hfill\mbox{}\\ 
Note that the $p$-adic Newton polygon on the left has exactly $2$ lower 
edges (with horizontal lengths $2$ and $3$), while the Archimedean Newton 
polygon on the right has exactly $3$ lower edges (with horizontal lengths 
$1$, $1$, and $3$). \dia}  
\end{ex}  

\begin{thm} \label{thm:newt} {\em 
Following the notation above, the number\footnote{counting multiplicity} 
of roots of $f$ in $\C_p$ of valuation $v$ is exactly the horizontal length of 
the face of $\newt_p(f)$ with inner normal $(v,1)$. Furthermore, if $\anewt(f)$ 
has a lower edge $E$ with slope $v$, and no other lower edges with slope 
in the open interval $(v-\log 3,v+\log 3)$, then the number$^3$ of roots 
$\zeta\!\in\!\C$ of $f$ with $\log|\zeta|\!\in\!(v-\log 3,v+\log 3)$ 
is exactly $\lambda(E)$. \qed }  
\end{thm} 

\noindent 
The first portion of Theorem \ref{thm:newt} goes back to early 
20th century work of Hensel, while the second portion is an 
immediate consequence of \cite[Thm.\ 1.5]{aknr} (with an  
important precursor in \cite{ostrowskiarch}). 

We will also use the following version of Hensel's famous 
criterion for the rapid convergence of Newton's method over $\C_p$: 
\begin{hensel} 
{\em (See, e.g., \cite[Thm.\ 4.1 \& Inequality (5.7)]{conrad}.)  
Suppose $p$ is prime, $f\in \Z[x]$, $j\!\geq\!1$,   
$\zeta\!\in\!\Z_p$, $\ell\!=\!\ord_p f'(\zeta)\!<\!\infty$, and 
$f(\zeta)\equiv 0\mod p^{2\ell+j}$. 
Let $\zeta'\!:=\!\zeta-\frac{f(\zeta)}{f'(\zeta)}$. 
Then $f(\zeta')\!=\!0$ mod $p^{2\ell+2j}$, $\ord_p f'(\zeta')\!=\!\ell$, 
and $\zeta\!=\!\zeta'$ mod $p^{\ell+j}$. \qed }  
\end{hensel}

\subsection{Separating Roots of Binomials} 
\label{sec:bisep}
When $f\!\in\!\Z[x]$ is a binomial, all of its roots in $\C$ are 
multiples of roots of unity that are evenly spaced on a circle. 
The same turns out to be true over $\C_p$, but the root spacing 
then depends more subtly on $p$ and less on the degree. For convenience, 
we will sometimes write $|\cdot|_\infty$ instead of $|\cdot|$ 
for the standard norm on $\C$. It will be convenient to equivalently rephrase 
lower bounds on distances between roots $|\zeta_1-\zeta_2|_p$ (which always 
tend to $0$ as $H\longrightarrow\infty$ in our setting)  
as {\em upper bounds} on $|\log|\zeta_1-\zeta_2|_p|$. 
In summary, we have the following: 
\begin{prop}{\em \label{prop:bi} 
Suppose $f(x)\!:=\!c_1+c_2x^d\!\in\!\Z[x]$, $c_1c_2\!\neq\!0$, and 
$|c_1|,|c_2|\!\leq\!H$. Also let $p\!\in\!\{\infty,2,3,5,\ldots\}$ 
and let $\overline{K}_p$ denote $\C$ or $\C_p$, according as $p\!=\!\infty$ 
or $p$ is prime. Then for any distinct roots
$\zeta_1,\zeta_2\!\in\!\overline{K}_p$ of $f$, we have that  
$|\log|\zeta_1-\zeta_2|_p|$ is bounded from above by: \\  
\mbox{}\hfill$\begin{cases} 
\log(d)+\frac{1}{d}\log H ; \text{ for } 
p\!=\!\infty \text{ and } d\!\geq\!2, \\ 
\mbox{}\hspace{1.55cm}\frac{1}{d}\log H; \text{  for } 
d\!>\!p^{\ord_p d}, \text{ and} \\ 
\mbox{}\hspace{.35cm}\frac{\log p}{p-1}+\frac{1}{d}\log H; 
\text{ for } d\!=\!p^{\ord_p d}\!\geq\!p. \end{cases}$\hfill\mbox{}}  
\end{prop} 

\noindent 
{\bf Proof:} Please see the Appendix, Section \ref{sub:propbi}. \qed 

\smallskip 
\noindent 
It is interesting that if one fixes $p$ and $H$, and lets 
$d\longrightarrow \infty$, then the minimal root distance tends to $0$ for the 
Archimedean case ($p\!=\!\infty$), but {\em is never less than}  
$\frac{1}{Hp^{1/(p-1)}}$ 
for the non-Archimedean case ($p$ prime).  

\subsection{Counting Roots of Binomials Over $\Q^*_p$} 
\label{sub:binocountqp}  
For any ring $R$ we let $R^*$ denote the 
multiplicatively invertible elements of $R$.  
Counting roots of binomials 
over $\Q_p$ is more involved than counting their roots over $\R$, but is 
still quite efficiently doable. 
\begin{lem} {\em \label{lem:binoqp} 
Suppose $p$ is an odd prime and $f(x)\!:=\!c_1+c_2x^d\!\in\!\Z[x]$ with 
$|c_1|,|c_2|\!\leq\!H$, $c_1c_2\!\neq\!0$, and $\ell\!:=\!\ord_p d$. 
Then the number of roots of $f$ in $\Q_p$ is either $0$ or 
$\gcd(d,p-1)$. In particular, $f$ has roots in $\Q_p$ if and only if 
{\em both} of the following conditions hold:\\ 
\mbox{}\hfill 
(1) $d|\ord_p(c_1/c_2)$ \ and \ (2) $\left(-\frac{c_1}{c_2}p^{\ord_p(c_2/c_1)}
\right)^{p^{\ell}(p-1)/\gcd(d,p-1)}\!=\!1$ {\em mod} 
$p^{2\ell+1}$. \hfill \qed }  
\end{lem} 

\noindent 
Lemma \ref{lem:binoqp} is classical and follows from basic group theory 
(the fact that the multiplicative group $(\Z/(p^k))^*$ is cyclic, of order 
$p^{k-1}(p-1)$, for $p$ odd) and Hensel's Lemma. The case $p\!=\!2$ is 
slightly more involved and is stated in the Appendix, Section \ref{sub:2adic}. 
\begin{cor} {\em \label{cor:binomod} 
Following the notation and assumptions of Lemma \ref{lem:binoqp}, 
one can count exactly the number of roots of $f$ in $\Q_p$ in 
time $O\!\left([\log(dpH)\log\log(dpH)]^2\right)$.  
Furthermore, for any root $\zeta\!\in\!\Q^*_p$ there is an 
$x_0\!\in\!\Z\left/\left( p^{2\ell+1}\right) 
\right.$ that is a root of the mod $p^{2\ell+1}$ reduction of\linebreak  
$\frac{c_1}{p^{\ord_pc_1}}+ \frac{c_2}{p^{\ord_pc_2}}x^d$,
and with $z_0\!:=\!p^{\ord_p(c_2/c_1)/d}x_0\!\in\!\Q$  
an approximate root of $f$ with associated true root $\zeta$. In 
particular, the logarithmic height\footnote{The logarithmic height of 
a rational number $a/b$ with $\gcd(a,b)\!=\!1$ is simply 
$\log\max\{|a|,|b|\}$ (and 
we declare the logarithmic height of $0$ to be $0$).} 
of $z_0$ is $O\!\left(\log\left(pH^{1/d} \right)\right)$. }  
\end{cor}  

\noindent 
{\bf Proof:} Please see the Appendix, Sections \ref{sub:oddbinomod} 
and \ref{sub:2binomod}. \qed

\subsection{Trees and Roots in $\Z/(p^k)$ and $\Z_p$}   
\label{sub:trees} 
The $p$-adic analogue of bisecting an isolating interval containing a real 
root is to approximate the next base-$p$ digit of an approximate root in 
$\Q_p$. Shifting from bisecting intervals to extracting 
digits is crucial since $\Q_p$ is not an\linebreak 
\scalebox{.96}[1]{ordered field. We will write $f'$ for 
the derivative of $f$ and $f^{(i)}$ for the $i$th order derivative of $f$.}  
\begin{dfn} {\em \label{dfn:crazytree} \cite{krrz19} 
For any $f\in \Z[x]$ let $\tilde{f}$ denote the 
mod $p$ reduction of $f$. A root $\zeta_0\!\in\!\F_p$ 
of $\tf$ is {\em degenerate} if and only if $\tf'(\zeta_0)\!=\!0$ mod $p$. 
For any degenerate root $\zeta_0$ of $\tilde{f}$ (represented as an 
element of $\{0,\ldots,p-1\}$), 
we then define $s(f,\zeta_0):=\min_{i\geq 0}\{ i
+\ord_p\frac{f^{(i)}(\zeta_0)}{i!}\}$.   
Fixing $k\in \N$, for $i\geq 1$, let us inductively define a set 
$T_{p,k}(f)$ of pairs $(f_{i-1,\mu}, k_{i-1,\mu})$   
$\in \Z[x]\times \N$: We set $(f_{0,0}, k_{0,0}) := (f,k)$. Then for any 
$i\geq 1$ with $(f_{i-1,\mu}, k_{i-1,\mu})$ 
$\in T_{p,k}(f)$, and any degenerate root $\zeta_{i-1}\!\in\!\{0,\ldots,p-1\}$  
of $\tilde{f}_{i-1,\mu}$ with
$s_{i-1}:= s(f_{i-1,\mu},\zeta_{i-1})\in \{2,\ldots,k_{i-1,\mu}-1\}$, 
we define $\zeta:= \mu + \zeta_{i-1}p^{i-1},k_{i,\zeta}:=k_{i-1,\mu} 
- s_{i-1}$, $f_{i,\zeta}(x) := p^{-s(f_{i-1,\mu},\zeta_{i-1})} 
f_{i-1,\mu}(\zeta_{i-1} + px)  \mod p^{k_{i,\zeta}}$, and 
then include $(f_{i,\zeta},k_{i,\zeta})$ in $T_{p,k}(f)$. \dia}  
\end{dfn}

\begin{ex} \label{ex:tri} 
{\em If $f(x)\!=\!x^{10}-10x+738$ and $p\!=\!3$ then $\tf(x)\!=\!x(x-1)^9$ mod 
$3$, $1$ is a degenerate root of $\tf$, and one can check that $s(f,1)\!=\!4$ 
(no greater than the multiplicity of the factor $x-1$ in $\tf$). In particular, 
$f_{1,1}$ has degree $10$ (and $10$ monomial terms) but 
$\tf_{1,1}\!=\!x^3+2x^2$. \dia } 
\end{ex} 

The collection of pairs $(f_{i,\zeta},k_{i,\zeta})$ admits a tree structure 
that will give us a way to extend Hensel lifting to degenerate roots. 
\begin{dfn}
\label{dfn:tree} {\em \cite{krrz19} 
Let us identify the elements of $T_{p,k}(f)$ with nodes of a labelled, rooted, 
directed tree $\cT_{p,k}(f)$ defined inductively as follows\footnote{This 
definition differs slightly from the original in \cite{krrz19}. 
}: 
\begin{itemize}
\item[(i)]{We set $f_{0,0}\!:=\!f$, $k_{0,0}\!:=\!k$, and let
$(f_{0,0},k_{0,0})$ be the label of the root node of
$\cT_{p,k}(f)$.}  
\item[(ii)]{\scalebox{.94}[1]{The non-root nodes of $\cT_{p,k}(f)$ are 
uniquely labelled by each $(f_{i,\zeta},k_{i,\zeta})$ 
$\in T_{p,k}(f)$ with $i\!\geq\!1$.}}
\item[(iii)]{There is an edge from node $(f_{i-1,\mu},k_{i-1,\mu})$ to
node $(f_{i,\zeta},k_{i,\zeta})$ if and only if 
there is a degenerate root $\zeta_{i-1}\!\in\!\{0,\ldots,
p-1\}$ of $\tf_{i-1,\mu}$ with $s(f_{i-1,\mu},\zeta_{i-1})
\!\in\!\{2,\ldots,k_{i-1,\mu}-1\}$ and 
$\zeta\!=\!\mu+\zeta_{i-1}p^{i-1}\!\in\!\Z/(p^i)$. \dia}
\end{itemize} } 
\end{dfn}

\noindent 
We call each $f_{i,\zeta}$ with $(f_{i,\zeta},k_{i,\zeta})\!\in\!T_{p,k}(f)$ a 
{\em nodal polynomial} of $\cT_{p,k}(f)$. 
It is in fact possible to list all the roots of $f$ in $\Z/(p^k)$ from the 
data contained in $\cT_{p,k}(f)$ \cite{krrz19,DMS19}. We 
will instead use $\cT_{p,k}(f)$, with $k$ chosen via our root separation 
bounds, to efficiently {\em count} the roots of $f$ in $\Z_p$, and then in 
$\Q_p$ by rescaling.  

\begin{ex} 
{\em $\cT_{p,k}(x^2)$ is a chain of length $\floor{\frac{k-1}{2}}$ for any 
$p,k$. \dia } 
\end{ex} 
\begin{ex} 
{\em Let $f(x)\!=\!1-x^{397}$. Then $\cT_{17,k}(f)$, for any $k\!\geq\!1$,  
consists of a single node, labelled $(1-x^{397},k)$, since $\tilde{f}$ has 
no degenerate roots in $\F_{17}$. In particular, $f$ has $1$ as 
its only root in $\Q_{17}$. 
\dia}   
\end{ex}            
\begin{ex}         
{\em Let $f(x)\!=\!1-x^{340}$. Then, when 
$k\!\in\!\{1,2\}$, the tree $\cT_{17,k}(f)$ consists 
of a single root node, labelled $(1-x^{340},k)$. 
However, when $k\!\geq\!3$, the tree $\cT_{17,k}(f)$ has depth $1$, 
and consists of the aforementioned root node {\em and} exactly $4$ child 
nodes, labelled $(f_{1,\zeta_0},k-2)$ where the 
$\tf_{1,\zeta_0}$ are, respectively, $14x$, $12x+10$, 
$5x+15$, and $3x+3$. Note that $\tilde{f}$ has 
exactly $4$ roots $\zeta_0\!\in\!\F_{17}$ ($1$, $4$, $13$, and $16$), 
each of which is degenerate, and the roots $\zeta_1\!\in\!\F_{17}$ of the 
$\tf_{1,\zeta_0}$ encode the ``next'' base-$17$ digits ($0$, $2$, 
$14$, and $16$) of the roots of $f$ 
in $\Z/(17^2)$. In particular, the roots of $f$ in $\Q_{17}$ are 
$1+0\cdot 17+\cdots$, $4+2\cdot 17+\cdots$,
$13+14\cdot 17+\cdots$, 
and $16+16\cdot 17+\cdots$ and are all {\em non}-degenerate. 
\dia}  
\end{ex} 

Nodal polynomials thus encode individual base-$p$ digits of roots of 
$f$ in $\Z_p$. Their degree also decays in a manner depending on root 
multiplicity\footnote{Over any field $K$, we define the {\em multiplicity  
of a root $\zeta\!\in\!K$} of $f\!\in\!K[x]$ as the greatest $m$ with  
$(x-\zeta)^m|f$ in $K[x]$.} 
as follows: 

\begin{lem} {\em \label{lem:nodal} 
\cite[Lem.\ 2.2 \& 3.6]{krrz19} Following the notation of Definition 
\ref{dfn:tree}, suppose $i\!\geq\!1$, $\zeta_{i-1}$ has multiplicity $m$ 
over $\F_p$, and
$(f_{i,\zeta},k_{i,\zeta})\!\in\!T_{p,k}(f)$. Then $\cT_{p,k}(f)$ has 
depth $\!\leq\!\floor{(k-1)/2}$\linebreak 
\scalebox{.9}[1]{and 
$\deg \tf_{i,\zeta}\!\leq\!s(f_{i-1,\mu},\zeta_{i-1})\!\leq\!\min\{k_{i-1,
\mu}-1,m\}$. Also, $f_{i,\zeta}(x)=
p^{-s}f(\zeta_0+\zeta_1p+\cdots +\zeta_{i-1}p^{i-1}+p^ix)$}\linebreak  
where $s\!:=\!\sum^{i-1}_{j=0}
s(f_{j,\zeta_0+\cdots+\zeta_{j-1} p^{j-1}},\zeta_j)$. \qed} 
\end{lem} 

Let $n_p(f)$ denote the number of non-degenerate roots in $\F_p$ of the mod $p$ 
reduction of $f$. 
\begin{lem} {\em \label{lem:ulift} 
If $f\!\in\!\Z[x]$, $D$ is the maximum of $\ord_p(\zeta_1-\zeta_2)$ over all 
$\zeta_1,\zeta_2\!\in\!\Z_p$ with $f(\zeta_1)\!=\!f(\zeta_2)\!=\!0\!\neq\!
\zeta_1-\zeta_2$, and $k\!\geq\!1+D$, then $f$ has 
exactly $\sum_{(g,j)\in T_{p,k}(f)} n_p(g)$ non-degenerate 
roots in $\Z_p$. }  
\end{lem}

\noindent 
{\bf Proof:} 
By Lemma \ref{lem:nodal}, 
$f_{i,\zeta}(x)\!=\!\frac{1}{p^s}f(\zeta_0+\cdots+\zeta_{i-1}p^{i-1}
+p^ix)$. So by Hensel's Lemma, any root $\zeta_i\!\in\!\F_p$ of 
$f_{i,\zeta}$ lifts to a unique non-degenerate root $\zeta_i+p\zeta_{i+1}
+\cdots\!\in\!\Z_p$ of $f_{i,\zeta}$. In other words, 
we obtain $\zeta_0+\zeta_1p+\cdots\!\in\!\Z_p$ as a root of $f$.  
Any sequence $(\zeta_0,\ldots,\zeta_{i-1})\!\in\!\F^i_p$ defined 
by a nodal polynomial $f_{i,\zeta}$ thus determines a unique root in 
$\Z_p$ of $f$, and we thus see that $\sum_{(g,j)\in T_{p,k}(f)} n_p(g)$ is a 
lower bound on the number of non-degenerate roots of $f$ in $\Z_p$. 

To see that we obtain {\em all} non-degenerate roots of $f$ in $\Z_p$ 
this way, note that the mod $p^k$ reduction of any root 
of $f$ in $\Z_p$ is a root of the mod $p^k$ reduction of $f$ in 
$\Z/(p^k)$. By the definition of $k$, the resulting map is an injection 
since distinct roots in $\Z_p$ must differ somewhere within their $1+D$ 
most significant digits. \qed

\subsection{Trees and Extracting Digits of Radicals} 
We prove the following useful lemma in 
Remark \ref{rem:low} of Section \ref{sec:trinosolqp}:  
\begin{lem} {\em \label{lem:binodepth} 
Suppose $f(x)\!=\!c_1+c_2x^d\!\in\!\Z[x]$ with $c_1c_2\!\neq\!0$ mod $p$ 
and $\ell\!:=\!\ord_p d$. Then every {\em non}-root nodal polynomial 
$f_{i,\zeta}$ of $\cT_{p,k}(f)$ satisfies $\deg \tf_{i,\zeta}\!\leq\!2$ or 
$\deg \tf_{i,\zeta}\!\leq\!1$, according as $p\!=\!2$ or $p\!\geq\!3$. 
In particular, $f(\zeta_0)\!=\!0$ mod $p$ for some 
$\zeta_0\!\in\!\{0,\ldots,p-1\} \Longrightarrow s(f,\zeta_0)\!\leq\!\ell+1$.} 
\end{lem} 

With our tree-based encoding of $p$-adic roots in place, we can now 
prove that it is easy to find approximate roots in $\Q_p$ for binomials 
when $p$ is fixed.  
\begin{thm} {\em
\label{thm:binoqp}
Suppose $f\!\in\!\Z[x]$ is a binomial
of degree $d$ with coefficients of absolute value at most $H$, $f(0)\!\neq\!0$, 
$\gamma\!=\!\gcd(d,\max\{2,p-1\})$, and $\{\zeta_1,\ldots,\zeta_\gamma\}$ is 
the set of roots of $f$ in $\Q_p$. Then in time 
$O\!\left(p\log(dp)\log\log(dp)+[\log(dpH)\log\log(dpH)]^2\right)$, 
we can find, for each $j\!\in\!\{1,\ldots,\gamma\}$, a $z^{(j)}_0\!\in\!\Q$ of 
logarithmic height $O\!\left(\log\left(dH^{1/d}\right)\right)$ 
that is an approximate root of $f$ with associated true root $\zeta_j$. } 
\end{thm}

An algorithm that proves Theorem \ref{thm:binoqp} when $p$ is odd is outlined 
below. Please see the Appendix, Section \ref{sub:binoq2} for the case 
$p\!=\!2$.\\  
\mbox{}\scalebox{.96}[.95]{\fbox{\mbox{}\hspace{.1cm}\vbox{
\begin{algor} {\em 
\label{algor:binoqp}
{\bf (Solving Binomial Equations Over $\pmb{\Q^*_p}$)}  
\mbox{}\\
{\bf Input.} An odd prime $p$ and 
$c_1,c_2,d\!\in\!\Z\setminus\{0\}$ with $|c_i|\!\leq\!H$ for all $i$. \\ 
{\bf Output.} A true declaration that $f(x)\!:=\!c_1+c_2x^d$ has 
no roots in $\Q_p$, or $z_1,$ 
$\ldots,z_\gamma\!\in\!\Q$ with\\ 
\mbox{}\hspace{1.8cm}\scalebox{.98}[1]{logarithmic 
height $O\!\left(\log\left(dH^{1/d}\right)\right)$ such that $\gamma\!=\!
\gcd(d,p-1)$, $z_j$ is an approximate}\\ 
\mbox{}\hspace{1.8cm}\scalebox{.96}[1]{root of $f$ with associated true root 
$\zeta_j\!\in\!\Q_p$ for all $j$, and the $\zeta_j$ are pair-wise distinct.}\\
{\bf Description.} \\ 
1: \scalebox{.95}[1]{If $\ord_p c_1\!\neq\!\ord_p c_2$ mod $d$ then say 
{\tt ``No roots in $\Q_p$!''} and {\tt STOP}.} \\
2: Let $\ell\!:=\!\ord_p d$ and replace $f$ with $f(x)\!:=\!c'_1+c'_2x^d$ 
where $c'_i\!:=\!\frac{c_i}{p^{\ord_p c_i}}$ for all $i$. \\
3: \scalebox{.93}[1]{If $\left(-\frac{c'_1}{c'_2}
\right)^{p^{\ell}(p-1)/\gamma}\!\!\!\!\!\!\!\!\!\neq\!1$ 
mod $p^{2\ell+1}$ then say {\tt ``No roots in $\Q_p$!''} and {\tt STOP}.}\\ 
4: \scalebox{.96}[1]{Let $\delta\!:=\!1$. If $d\!\leq\!-1$ then 
set $\delta\!:=\!-1$ and respectively replace $d$ by $|d|$ and $f(x)$ by 
$x^df(1/x)$.}\\ 
5: \scalebox{.91}[1]{Let $g$ be any generator for $\F^*_p$, 
$r\!:=\!(d/\gamma)^{-1}$ mod $(p-1)$, 
$c'\!:=\!(-c'_1/c'_2)^r$ mod $p$, and $\tilde{h}(x)\!:=\!x^{\gamma}-c'$.}\\ 
6: Find a root $x_1\!\in\!\left\{g^0,\ldots,g^{\frac{p-1}{\gamma}-1}\right\}$ 
of $\tilde{h}$ via brute-force search.\\ 
7: For all $j\!\in\!\{2,\ldots,\gamma\}$ let 
$x_j\!:=\!x_{j-1}g^{(p-1)/\gamma}$ mod $p$. \\ 
8: \scalebox{.87}[1]{If $\ell\!\geq\!1$ then, for each $j\!\in\!\{1,\ldots,
\gamma\}$, replace $x_j$ by $x_j -\frac{f(x_j)/p^\ell}
{f'(x_j)/p^\ell}\!\in\!\Z/(p^2)$.}\\ 
9: Output $\left\{(x_1 p^{\ord_p(c_1/c_2)/d})^\delta,\ldots,(x_\gamma 
p^{\ord_p(c_1/c_2)/d})^\delta\right\}$. } 
\end{algor}} 
}}

\begin{rem} {\em 
We will see below that, upon 
rescaling all the roots and approximate roots to the $p$-adic unit circle,  
each approximate root is within $1/p$ or $1/p^2$ of a unique true root, 
according as $\ell$ is $0$ or not. Also, Step 6 is stated for simplicity 
rather than practicality, and\linebreak 
\scalebox{.98}[1]{can be sped up considerably if one 
one avails to randomization: See, e.g., \cite[Ch.\ 7, Sec.\ 3]{bs}. \dia}} 
\end{rem} 

\noindent
{\bf Proof of Theorem \ref{thm:binoqp}:} 
It clearly suffices to prove the correctness of Algorithm \ref{algor:binoqp}, 
and then analyze its complexity. In particular, we assume $p$ is odd 
in this proof. (Please see the Appendix, Section \ref{sub:binoq2} for the case 
$p\!=\!2$.) 

\medskip 
\noindent 
{\bf Correctness:} Theorem \ref{thm:newt} implies that Step 1 merely 
checks whether the valuations of the roots of $f$ in $\C^*_p$ in fact lie in 
$\Z$, which is necessary for $f$ to have roots in $\Q^*_p$. 

Steps 2 and 4 merely allow us to reduce our search for approximate roots to 
$(\Z/(p^{2\ell+1}))^*$ and assume positive degree $d$. 

\scalebox{.96}[1]{Lemma \ref{lem:binoqp} implies that 
Step 3 merely check that the coset of roots of $f$ in $\C^*_p$ 
intersects $\Z^*_p$.}  

Step 5 is merely the application of an automorphism of $\F^*_p$  
(that preserves the roots of $\tf$ in $\F^*_p$) that 
enables us to work with a binomial of degree $\gamma$ (possibly 
much smaller than both $p-1$ and $d$). 

Steps 6--7 then clearly find the correct coset of $\F^*_p$ that makes 
$f$ vanish mod $p$. In particular, by Hensel's Lemma, Step 9 clearly gives 
the correct output if $\ell\!=\!0$. (Recall that we have replaced 
each coefficient $c_i$ of $f$ with $c'_i$.) 

If $\ell\!\geq\!1$ then let $\zeta_0$ be any $x_j$ from Step 8. We then 
have $\deg \tilde{f}_{1,\zeta_0}\!\leq\!1$ thanks 
to Lemma \ref{lem:binodepth}. 
Furthermore, Definition \ref{dfn:crazytree} tells us that the unique root 
$\zeta_1\!\in\!\F_p$ of $\tilde{f}_{1,\zeta_0}$ is exactly the next 
base-$p$ digit of a unique root $\zeta\!\in\!\Z_p$ of $f$ with 
$\zeta\!=\!\zeta_0$. 
Also, $\deg \tilde{f}_{1,\zeta_0}$ must be $1$ (for otherwise $\tf$ would 
not vanish on its coset of roots in $\F^*_p$) and $s(f,\zeta_0)\!\geq\!2$ 
since $\ell\!\geq\!1$ forces $\zeta_0$ to be a degenerate root of $\tf$. 
Lemma \ref{lem:nodal} then tells us that Hensel's Lemma 
--- applied to $f_{1,\zeta_0}(x)\!=\!p^{-s(f,\zeta_0)}f(\zeta_0+px)$ and 
start point $\zeta_1\!\in\!\Z/(p)$ --- implies that $\zeta_0+\zeta_1p$ 
is an approximate root of $f$ with associated true root 
$\zeta\!\in\!\Z_p$. So Step 8 in fact refines $x_1$ to the mod $p^2$ 
quantity $\zeta_0+\zeta_1p$, and thus Steps 7--9 indeed give us suitable 
approximants in $\Q$ to all the roots of $f$ in $\Q_p$. So our algorithm 
is correct.  

Note also that the outputs, being integers in $\{0,\ldots,p^2-1\}$ 
rescaled by a factor of $p^{\ord_p(c_1/c_2)/d}$ (or possibly the reciprocals 
of such quantities), clearly each have bit-length $O\!\left(\log(p)
+\frac{|\log(c_1/c_2)|}{d\log p}\log p\right)\!=\!
O\!\left(\log(p)+\frac{\log H}{d}\right)\!=\!O\!\left(\log\!\left(
pH^{1/d}\right)\right)$. \qed  

\medskip 
\noindent 
{\bf Complexity Analysis:} Via Corollary \ref{cor:binomod}, \cite{generator}, 
and some additional elementary bit complexity estimates for modular 
arithmetic \cite{vzgbook}, it is clear that, {\em save for Steps 6--9}, 
Algorithm \ref{algor:binoqp} has complexity $O(p^{1/4}\log(p)\log\log(p)  
 +[\log(dpH)\log\log(dpH)]^2)$, 
provided we use Harvey and van der Hoeven's recent fast multiplication 
algorithm \cite{harveymult}. Steps 6--7 (whose complexity dominates the 
complexity of Steps 6--9), involve  $\frac{p-1}{\gamma}-1$ multiplications in 
$\F_p$ and $\gamma-1$ multiplications in $\Z/(p^{2\ell+1})$. 
This takes time no worse than $O(p\log(dp)\log\log(dp))$, so we are done. \qed 

\section{Proving Theorem \ref{thm:tri}}  
\label{sec:trisepqp} 
Let us first recall the following version of {\em Yu's Theorem}:
\begin{thm} \label{thm:yu}
\cite[Cor.\ 1]{yu94} {\em 
Suppose $p$ is any prime, $n\!\geq\!2$, $\alpha_1,\ldots,\alpha_n\!\in\!\Q$  
with $\alpha_i = r_i/s_i$ a reduced fraction for 
each $i$, and $b_1,\ldots,b_n\!\in\!\Z$ are not all zero. Then 
$\alpha_1^{b_1}\cdots\alpha_n^{b_n} \neq 1$ implies that
$\alpha_1^{b_1}\cdots\alpha_n^{b_n} - 1$ has $p$-adic valuation bounded 
from above by\\  
\mbox{}\hspace{.25cm}$11145\left( \frac{24(n+1)^2}{\log p}\right)^{n+2}(p-1)
\left(\prod^n_{i=1}\log A_i\right)\log(4B)
\times \max\left\{\log(2^{12}\cdot 3n(n+1)
\log A_n), \frac{\log p}{n} \right\}$,\\   
\scalebox{.96}[1]{where $B\!:=\!\max\{|b_1|,\ldots,|b_n|,3\}$, and 
$A_1,\ldots,A_n$ are any real numbers
such that $A_1\leq \cdots\leq A_n$}\linebreak 
and, for each $j$, $A_j \geq \max\{|r_j|,|s_j|,p\}$. \qed } 
\end{thm}

To prove that two distinct roots $\zeta_1,\zeta_2\!\in\!\C_p$ of a trinomial 
$f$ can not be too close, we will build a special point $m\!\in\!\C$ with four 
special properties: (i) $f'(m)\!=\!0$, (ii) $|f(m)|_p$ is not too small,
(iii) $|\zeta_1-\zeta_2|_p\!\geq\!p^{-1/(p-1)}|\zeta_1-m|_p$, and (iv) 
$|\zeta_1-m|$ is not too small. So let us quantify this approach toward  
proving Theorem \ref{thm:tri}.  
\begin{prop} {\em \label{prop:1}
Let $f(x)\!=\!c_1+c_2x^{a_2} +c_3x^{a_3}\!\in\!\Z[x]$ with
$a_3\!>\!a_2\!\geq\!1$, $c_1\!\neq\!0$, and suppose $m\!\in\!\C_p$ is a root 
of $f'$. Then $m^{a_3-a_2}\!=\!-\frac{a_2 c_2}{a_3 c_3}$ and 
$f(m)\!=\!c_1+c_2m^{a_2}\left(1-\frac{a_2}{a_3}\right)$. \qed } 
\end{prop}
\begin{lem}{\em 
\label{lem:1} 
Following the notation above, assume further that 
$f$ is square-free. Then\\  
$\norm{f(m)}\!\geq\!\exp\left(-O\!\left(\frac{p}{\log^2 p}\log(d)\log^2(dH+p) 
\log\log(dH+p)\right)\right)$. }  
\end{lem}

\noindent 
{\bf Proof:} First note that if $f$ is square-free then $f$ has no 
repeated factors, and thus no degenerate roots in $\C_p$. So $f(m)\!\neq\!0$. 
By Proposition \ref{prop:1} we then obtain $\ord_p f(m)$\linebreak   
\scalebox{.85}[1]{$= \ord_p(c_1+c_2m^{a_2}(1-a_2/a_3))
=\ord_p(c_1) + \ord_p(-1) + 
\ord_p\left(\frac{-c_2(a_3-a_2)}{c_1 a_3}
\left(-\frac{a_2 c_2}{a_3 c_3}\right)^{a_2/(a_3-a_2)}-1\right)$.} 
 \hfill (1) \label{eq:sum} 

Clearly, $\ord_p c_1 \leq \frac{\log H}{\log p}$ and $\ord_p(-1)\!=\!0$. 
To bound the third summand above, let $T\!:=\!\frac{-c_2(a_3-a_2)}{c_1 a_3}
\left(-\frac{a_2 c_2}{a_3 c_3}\right)^{a_2/(a_3-a_2)}$ 
and observe that $T^{a_3-a_2}-1\!=\!\prod^{a_3-a_2}_{j=1}(T-\zeta^j)$  
for $\zeta\!\in\!\C_p$ a primitive $(a_3-a_2)$-th root 
of unity. In particular, $T^{a_3-a_2}\!\neq\!1$ since $f(m\zeta^j)\!\neq\!0$ 
for all $j\!\in\!\{1,\ldots,a_3-a_2\}$, thanks to Proposition \ref{prop:1} 
and $f$ not having any degenerate roots. So then 
$M\!:=\!\ord_p(T^{a_3-a_2}-1) = \sum_{j=1}^{a_3-a_2}\ord_p(T-\zeta^j)
<\infty$, with the $(a_3-a_2)$-th term of the sum exactly 
$\ord_p(T-\zeta^{a_3-a_2})\!=\!\ord_p(T-1)$, 
i.e., the third summand from (1). 

Suppose $\ord_p T\!<\!0$. Then for each $i\!\in\!\{1,\ldots,a_3-a_2\}$, the 
Ultrametric Inequality gives us $\ord_p(T-\zeta^j) = \ord_pT < 0$, since roots 
of unity always have $p$-adic valuation $0$. We must then have 
$\ord_p f(m) = \ord_p(c_1) + \ord_p(T-\zeta^{a_3-a_2}) 
< \frac{\log(dH)} {\log p}$ and we obtain our lemma.

On the other hand, should $\ord_pT \geq 0$, we get
$\ord_p(T-\zeta^j)\geq j\ord_p(\zeta) = 0$, for each $j$, by the Ultrametric 
Inequality. So $M\!\geq\!\ord_p(T-1)$ and we'll be done if we find a 
sufficiently good upper bound on $M$.

By luck, $M$ is boundable directly from Theorem \ref{thm:yu}, upon setting   
$n\!=\!2$, $\alpha_1\!=\!-c_2\frac{a_3-a_2}{c_1 a_3}$, 
$\alpha_2\!=\!-\frac{a_2 c_2}
{a_3 c_3}$, $b_1\!=\!a_3-a_2$, and $b_2\!=\!a_2$. In particular, 
we can use $A_i = \max\{dH,p\}$ for $i\!\in\!\{1,2\}$ and $B = \max\{d,3\}$, 
yielding $\log A_1,\log A_2, \log B = O(\max\{\log(dH),\log p\})$, so that 
$M\!\leq\!Cp\log^2\max\{dH,p\}\log(4\max\{d,3\})
\left. \times \max\left\{\log\left(18\cdot 2^{12}
\log\max\{dH,p\}\right),\frac{\log p}{2}\right\}\right/\log^4 p$\\
for $C\!=\!11145\cdot 216^4$. So then 
$M\!=\!O\!\left(\frac{p}{\log^4 p}\log^2(dH+p)(\log d)(\log\log(dH+p)+\log p)
\right)$\\  
\mbox{}\hspace{3cm}$=O\!\left(\frac{p}{\log^4 p}\log(d)\log^2(dH+p)
\log\log(dH+p)\log p \right)$\\  
\mbox{}\hspace{3cm}$=O\!\left(\frac{p}{\log^3 p}\log(d)\log^2(dH+p)
\log\log(dH+p) \right)$.\\  
In other words, the third summand from (1) is bounded from above by 
the last $O$-bound,\linebreak 
\scalebox{.955}[1]{and thus $\ord_p f(m)\!=\!O(M)$ since $\frac{\log H}
{\log p}\!=\!O(M)$. Since $|f(m)|_p\!=\!e^{-\log(p)\ord_p f(m)}$, 
we are done. \qed}  

\medskip 
The Ultrametric Inequality directly yields the following:  
\begin{prop}{\em \label{prop:dev_small} 
If $f\!\in\!\Z[x]$ and $r\!\in\!\C_p$ then 
$\norm{r}\leq 1\Longrightarrow \norm{f'(r)}\leq 1$. \qed } 
\end{prop}

\medskip 
\scalebox{.96}[1]{Below is a rescaled {\em $p$-adic} version of 
{\em Rolle's Theorem}, based on \cite[Sec.\ 2.4, Thm., Pg.\ 316]{Rob00}.}   
\begin{thm}\label{thm:rolle} {\em Let $f\in \mathbb{C}_p[x]$ have two distinct roots $\zeta_1, \zeta_2\in \mathbb{C}_p$ 
with $\norm{\zeta_1-\zeta_2} = cp^{1/(p-1)}$ for some $c>0$. 
Then $f'$ has a root $m\!\in\!\C_p$ with $|\zeta_1-m|_p,|\zeta_2-m|_p\!\leq\!c$. \qed } 
\end{thm}

We can now prove one of our main results. 

\medskip 
\noindent 
{\bf Proof of Theorem \ref{thm:tri}:} For convenience, let us abbreviate 
the stated $O$-bound by $O(M)$. Note that if one of the $\zeta_i$ is $0$ then 
the other root $\zeta_j$ is nonzero, with valuation satisfying 
$|\ord_p \zeta_j|\!\leq\!\frac{\log H}{\log p}$, thanks to Theorem 
\ref{thm:newt}. So then 
$|\log|\zeta_i-\zeta_j|_p|\!=\!|\log|\zeta_j|_p|\!=\!\left|\log 
e^{-\log(p)\ord_p \zeta_j}\right|
=\!|-\log(p)\ord_p\zeta_j|\!\leq\!\log H \!=\!O(M)$. 
So we may assume $\zeta_1\zeta_2\!\neq\!0\!\neq\!f(0)$. 

\medskip 
\noindent 
{\bf Case 1: (Both roots are small: $\pmb{\norm{\zeta_1}, 
\norm{\zeta_2}\leq 1}$.)}  \\
Suppose $\norm{\zeta_1-\zeta_2} > p^{-2/(p-1)}$. Then  
$\norm{\zeta_1-\zeta_2} > e^{-2\log(p)/(p-1)}$. Since 
$2\log(p)/(p-1)=O(M)$ we are done.   

Now assume that $\norm{\zeta_1-\zeta_2}\leq p^{-2/(p-1)}$. Then by Theorem  
\ref{thm:rolle} there is an $m\in \mathbb{C}_p$ such that $f'(m) = 0$ and 
$\norm{\zeta_i-m} \leq p^{1/(p-1)}\norm{\zeta_1-\zeta_2} 
\leq p^{-1/(p-1)}$ for all $i\!\in\!\{1,2\}$. 
Note that the Ultrametric Inequality implies that 
$|m|_p\!\leq\!p^{-1/(p-1)}$.

Since $f$ is square-free, Lemma \ref{lem:1} implies that 
$\norm{f(m)} \geq e^{-O(M)}$. Applying Theorem 
\ref{thm:rolle} 
to $g(x)\!:=\!f(x)-\frac{f(m)-f(\zeta_1)}{m-\zeta_1}x-\frac{mf(\zeta_1)
-\zeta_1f(m)}{m-\zeta_1}$ 
(which vanishes at $m$ and $\zeta_1$), we then see that there is a 
$\zeta\!\in\!\C_p$ with $\norm{\zeta-\zeta_1}\leq 1$ (and thus 
$|\zeta|_p\!\leq\!1$) such that $g'(\zeta)\!=\!0$, i.e., 
$f(m) = f(m) - f(\zeta_1) = f'(\zeta)(m-\zeta_1)$. As $f(m) \neq 0$ we get 
$f'(\zeta)\neq 0$ and $m\neq \zeta_1$. From Proposition \ref{prop:dev_small} 
we have $\norm{f'(\zeta)} \leq  1$, so then $\norm{m-\zeta_1} = 
\frac{\norm{f(m)}}{\norm{f'(\zeta)}} \geq e^{-O(M)}$. 
We thus get  $\norm{\zeta_1-\zeta_2} \geq p^{-1/(p-1)}\norm{m-\zeta_1}\geq 
e^{-O(M)-\frac{\log p}{p-1}}\!=\!e^{-O(M)}$. \qed 

\medskip 
\noindent
\scalebox{.97}[1]{{\bf Case 2: (Both roots are large: $\pmb{\norm{\zeta_1}, 
\norm{\zeta_2} > 1}$.)} Please see the Appendix, Section \ref{sub:bigroots}. 
\qed}   

\medskip  
\noindent
{\bf Case 3: (Only one root has norm $\pmb{>1}$.)}  \\ 
Without loss of generality, we may assume that $|\zeta_1|_p\!\leq\!1\!<\!
|\zeta_2|_p$. We then simply note that, as $\norm{\zeta_1}\neq \norm{\zeta_2}$, 
we have $\norm{\zeta_1-\zeta_2} = \max\left\{\norm{\zeta_1}, \norm{\zeta_2}
\right\}\!>\!1$  and we are done. \qed 

\section{Proving Theorem \ref{thm:tetra}} 
\label{sec:tetra}

\subsection{The Case of Prime $\pmb{p}$} 
Let $g(x) = p^{2h}f(x+p^{h-1}) = p^{2h}(x+p^{h-1})^d 
- p^{2h}\left(\frac{x+p^{h-1}}{p^h} 
- \frac{1}{p}\right)^2$  
$=p^{2h}(x+p^{h-1})^d - x^2$. 
Then $g$ has the same roots as $f_{d,p}$, save for a ``small'' shift by 
$p^{h-1}$. Rescaling, we get \scalebox{.95}[1]
{$G(x):=\frac{g(p^{(h-1)d/2+h}x)}{p^{(h-1)d+2h}}
=p^{-(h-1)d-2h}\left[ p^{2h}(p^{(h-1)d/2+h}x+p^{h-1})^d - 
p^{(h-1)d+2h}x^2\right]$}\linebreak 
$=\sum_{i=0}^{d}{d\choose i}p^{(h-1)(di/2-i)+ih}x^i - x^2=
1-x^2 \mod p^{d(h-1)/2+1}$, 
which is square-free for odd prime $p$. 
(The case of $p\!=\!2$ is in the Appendix.)  of the 
Hensel's Lemma then implies that there are 
roots $\zeta_1,\zeta_2 \in \Z_p$ of $G$ such that 
$\zeta_1 \equiv 1 \mod p^{d(h-1)/2+1}$ and $\zeta_2 \equiv -1 \mod 
p^{d(d-1)/2+1}$. 

So $\norm{\zeta_1} = \norm{\zeta_2} = 1$. For each 
$i\!\in\!\{1,2\}$, $y_i = p^{(h-1)d/2+h}\zeta_i$ is the corresponding root of 
$G$, and thus of $g$. Then $x_1 = y_1 + p^{h-1}$ and $x_2 = y_2 + p^{h-1}$ are 
two roots of $f$ in $\Z_p$ such that
$\norm{x_1-x_2} = \norm{(y_1+p^{h-1}) - (y_2+p^{h-1})} = \norm{y_1-y_2} 
\leq \max\left\{\norm{y_1}, \norm{y_2}\right\} = 
p^{-(h-1)d/2-h} = p^{-\Omega(dh)}.$ \qed  

\subsection{The Case $\pmb{p\!=\!\infty}$} 
Please see the Appendix, Section \ref{sub:infinity}. 

\section{Solving Trinomials over $\Q_p$} \label{sec:trinosolqp} 
Unlike the binomial case, the tree $\cT_{p,k}(f)$ can have high 
depth for large $k$ and $f$ an arbitrary trinomial. However, Lemma 
\ref{lem:low} below will show that 
the structure of $\cT_{p,k}(f)$ is still simple: Depth no greater than 
$\floor{(k-1)/2}$, and all but possibly one path in $\cT_{p,k}(f)$ having 
no more than $2$ vertices of out-degree $\geq\!2$. We will prove an upper 
bound on $k$ that is large enough to count all roots in $\Z_p$ (via Lemma 
\ref{lem:ulift}), but still small enough for us to approximate all these roots 
in time $p^{5+o(1)}\log^{16+o(1)}(dH)$.  

We begin with a central bound, derived via Theorem \ref{thm:yu}: 
\begin{thm} {\em \label{thm:air} 
If $f(x)\!=\!c_1+c_2x^{a_2}+c_3x^{a_3}\!\in\!\Z[x]$ is a trinomial 
of degree $d\!=\!a_3\!>\!a_2\!\geq\!1$, with 
coefficients of absolute value at most $H$, then 
$\!\!\!
\displaystyle{\sum\limits_{\zeta\in\Z_p \; : \; f(\zeta)=\ord_p \zeta= 0\neq 
f'(\zeta)}} 
\!\!\!\!\!\!\!\!
\!\!\!\!\!\!\!\!
\ord_p f'(\zeta)\!=\!O(p^2\log^8(dH))$. } 
\end{thm} 

\noindent 
{\bf Proof:} Please see the Appendix, Section \ref{sub:gcd}. \qed 

\begin{lem} {\em \label{lem:low} 
Following the notation and assumptions of Theorem \ref{thm:air}, every 
{\em non}-root  
nodal polynomial $f_{i,\zeta}$ of $\cT_{p,k}(f)$ with $\zeta\!\neq\!0$ 
mod $p$ satisfies $\deg \tf_{i,\zeta}\!\leq\!4$, $\deg 
\tf_{i,\zeta}\!\leq\!3$, or $\deg \tf_{i,\zeta}\!\leq\!2$, according as 
$p\!=\!2$, $p\!=\!3$, or $p\!\geq\!5$.}  
\end{lem} 
\begin{ex} {\em \label{ex:upperbound} Recalling Example \ref{ex:tri}, which 
had $f(x)\!=\!x^{10}-10x+738$, 
observe that $\cT_{3,7}(f)$ is a chain of length $2$. In particular, 
$\tf_{1,1}(x)\!=x^2(x-1)$, $0$ is a degenerate root of $\tf_{1,1}$, and 
$s(f_{1,1},0)\!=\!2$. We can then easily calculate that 
$\tf_{2,1}(x)\!=\!2(x-1)(x-2)$ mod $3$.

There are a total of $4$ non-degenerate roots in $\F_3$ for the nodal 
polynomials: $1$ for $\tf_{0,0}$, $1$ for $\tf_{1,1}$, and $2$ for $\tf_{2,1}$. 
These non-degenerate roots in $\F_3$ then lift to the following 
roots of $f$ in $\Z_3$: $0+O(3^1)$,
$1+1\cdot 3 + O(3^2)$, $1+0\cdot 3+ 1\cdot3^2+O(3^3)$, and $1+0\cdot 3 + 
2\cdot 3^2+O(3^3)$. 
A quick calculation via Maple's {\tt rootp} command tells us 
that these are all the $3$-adic rational roots of $f$. \dia} 
\end{ex}
\begin{ex} 
{\em One can check that for $f(x)\!:=\!x^{10}+11x^2-12$, the tree 
$\cT_{2,8}(f)$ is isomorphic to \raisebox{-.15cm}
{\epsfig{file=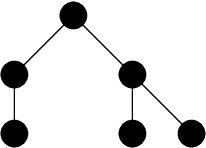,height=.5cm,clip=}}. 
In particular, this $f$ has exactly $6$ roots in $\Q^*_2$: 
$\tf_{2,2}\!=\!\tf_{2,1}\!=\!\tf_{2,3}\!=\!x^2+x$ and  
each of these (terminal) nodal polynomials has exactly $2$ 
non-degenerate roots in $\F_2$. Remembering the earlier digits encoded 
in $\cT_{2,8}(f)$, these $6$ roots then each lift to a unique root of $f$ in 
$\Z_2$. Note that $\tf_{1,1}(x)\!=\!x^4+x^2$ has degree $4$. \dia} 
\end{ex} 
\begin{ex}
{\em Composing Example \ref{ex:tri} with $x^2$, let us take 
$f(x)\!:=\!x^{20}-10x^2+738$. One then sees that the tree
$\cT_{3,7}(f)$ is isomorphic to \raisebox{-.15cm}
{\epsfig{file=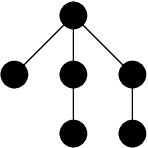,height=.5cm,clip=}}.
In particular, this $f$ has exactly $8$ roots in $\Q^*_3$, each arising 
as a Hensel lift of a non-degenerate root in $\F_3$ of some nodal polynomial: 
$\tf_{1,0}$, $\tf_{1,1}$, $\tf_{2,1}$, $\tf_{1,2}$, and $\tf_{2,8}$ 
respectively contribute $2$, $1$, $2$, $1$, and $2$ roots. 
Note that $\tf_{1,2}(x)\!=\!x^3+2x^2+x$ has degree $3$. \dia}
\end{ex}

To prove Lemma \ref{lem:low} we will need a powerful result of Lenstra  
\cite{len99} on the Newton polygons of shifted sparse polynomials. First, let 
us define $d_m(r)$ to be the least common multiple of all integers that 
can be written as the product of at most $m$ pairwise distinct positive 
integers that are at most $r$, and set $d_m(r)\!:=\!1$ if $mr\!=\!0$. 
\begin{thm} \cite[Sec.\ 3]{len99}
\label{thm:lenstra} 
{\em Suppose $f\!\in\!\Q[x]$ is a $t$-nomial, $g(x)\!=\!f(1+px)$, and 
$r$ is the largest nonnegative integer such that
$r-\ord_pd_{t-1}(r)\!\leq\!\underset{0\leq j\leq t-1}{\max}\{j-\ord_p(j!)\}$. 
Then any lower edge of 
$\newt_p(g)$ with inner normal $(v,1)$ with 
$v\!\geq\!1$ lies in the strip $[0,r]\times \R$. \qed}   
\end{thm} 

\noindent 
We point out that the vector of parameters $(t,r,v)$ from our statement 
above would be $(k+1,m,\nu(x-1))$ in the notation of \cite{len99}, and the 
parameter $r$ there is set to $1$ in our application here. 

\medskip 
\noindent 
{\bf Proof of Lemma \ref{lem:low}:} First note that replacing 
$x$ by $cx$, for any $c\!\in\!\{1,\ldots,p-1\}$, preserves the number of roots 
of $f$ in $\Z_p$ and (up to relabelling the $\zeta$ in the subscripts of the  
$f_{i,\zeta}$) the tree $\cT_{p,k}(f)$. So 
to study $\tf_{1,\zeta_0}$ with $\zeta_0\!\in\!\{1,\ldots,p-1\}$,  
it suffices to study $\tf_{1,1}$. 

Note that the lower hull of any Newton polygon can be identified with a 
piecewise linear convex function on an interval. In particular, 
$f_{1,1}(x)\!=\!p^{-s(f,1)}f(1+px)$ and thus the lower hull of 
$\newt_p(f_{1,1})$ can be identified with the sum of the lower hull of 
$\newt_p(f(1+x))$ and the function $x-s(f,1)$. Note also that by the 
definition of $\newt_p$, the minimal $y$-coordinate of a point 
of $\newt_p(f(1+px))$ is exactly $s(f,1)$. In particular, if 
$f(x)\!=\!c_1+c_2x^d$ with $c_1c_2\!\neq\!0$ and $f(\zeta_0)\!=\!0$ mod $p$ 
for some $\zeta_0\!\in\!\{1,\ldots,p-1\}$,  
then the definition of $s(f,\zeta_0)$ tells us that 
$s(f,\zeta_0)\!\leq\!1+\ord_p f'(\zeta_0)\!=\!1+\ord_p d\!=\!1+\ell$. 

Theorem \ref{thm:lenstra} then 
tells us that all lower edges of $\newt_p(f_{1,1})$ of non-positive slope 
lie in the strip $[0,r]\times \R$, where $r$ is the largest 
nonnegative integer such that 

\smallskip 
\noindent 
\mbox{}\hspace{.5cm} ($\star$) \hspace{4.5cm} 
$r-\ord_pd_2(r)\!\leq\!\eps_p$,

\smallskip 
\noindent 
where $\eps_2\!=\!1$ and $\eps_p\!=\!2$ for all $p\!\geq\!3$. In particular, 
the definition of $\newt_p(f_{1,1})$ tells us that $p$ divides the coefficient 
of $x^j$ in $f_{1,1}$ for all $j\!\geq\!r+1$ and thus 
$\deg \tf_{1,1}\!\leq\!r$. 

By Lemma \ref{lem:nodal}, all other non-root nodal 
polynomials $f_{i,\zeta}$ with $\zeta\!\neq\!0$ mod $p$ satisfy 
$\deg \tf_{i,\zeta}\!\leq\!\deg \tf_{1,1}$. So it 
suffices to prove that $r$ satisfies the stated bounds of our lemma. 
This is easily verified by first observing that $d_2(0)\!=\!d_2(1)\!=\!1$ 
and $d_2(2)\!=\!2$. So Inequality ($\star$) certainly holds for 
$r\!\in\!\{0,1,2\}$, regardless of $p$. Observing that $d_2(3)\!=\!6$ 
and $d_2(4)\!=\!24$, we then see that Inequality ($\star$) 
holds at $r\!=\!4$ (resp.\ $r\!=\!3$) when $p\!=\!2$ (resp.\ 
$p\!=\!3$). 

So it is enough to show that:\\   
\mbox{}\hspace{3.7cm}(i) $r-\ord_2 d_2(r)\!\geq\!2$ for $r\!\geq\!5$,\\  
\mbox{}\hspace{3.65cm}(ii) $r-\ord_3 d_2(r)\!\geq\!3$ for $r\!\geq\!4$, and \\ 
\mbox{}\hspace{3.54cm}(iii) $r-\ord_p d_2(r)\!\geq\!3$ for $r\!\geq\!3$ and 
$p\!\geq\!5$.\\ 
From \cite[Prop.\ 2.4]{len99}, we have $\ord_p d_2(r)\!\leq\!2
\frac{\log r}{\log p}$. Note that, for any fixed $p$, the 
quantity $r-2\frac{\log r}{\log p}$ is an increasing function of $r$ 
for $r\!\geq\!\frac{2}{\log p}$. Furthermore, 
$\ceil{7-2\frac{\log 7}{\log p}}\!\geq\!2$ 
for all $p\!\geq\!2$ and $\ceil{5-2\frac{\log 5}{\log p}}\!\geq\!3$
for all $p\!\geq\!3$. Noting that $d_2(5)\!=\!120$ and $d_2(6)\!=\!360$, it is 
then easily checked that (i)--(iii) all hold. \qed 
\begin{rem} 
\label{rem:low} 
{\em The proof of Lemma \ref{lem:binodepth} is simply the variation of 
the proof above where we replace Inequality ($\star$) by 
$r-\ord_pd_1(r)\!\leq\!1$, replace $d_2(r)$ with $d_1(r)$, and 
$\eps_p\!=\!1$ for {\em all} $p$. \dia } 
\end{rem}

\begin{lem} {\em \label{lem:complexity}
For any trinomial $f\!\in\!\Z[x]$ with $\tf(0)\!\neq\!0$ mod $p$, 
we can compute the mod $p$ reductions of 
all the nodal polynomials of $\cT_{p,k}(f)$ in time $p^{1+o(1)}k^{2+o(1)}$.}  
\end{lem}
\noindent 
{\bf Proof:} By Lemma \ref{lem:nodal}, 
$\cT_{p,k}(f)$ has depth $\leq\!\floor{\frac{k-1}{2}}$. 
By Lemma \ref{lem:low}, all non-root nodal polynomials have mod $p$ reduction 
of degree no greater than $4$. Thus, the root of $\cT_{p,k}(f)$ has at most 
$p-1$ children (since $\tf(0)\!\neq\!0$), and any node at depth $\geq\!1$ has 
no more than $2$ children (since a degree $4$ polynomial has at most $2$ 
degenerate roots). Lemma \ref{lem:nodal} also tells us that 
$\deg \tf_{i,\mu+\zeta_{i-1}p^{i-1}}$ is at most the multiplicity of 
$\zeta_{i-1}\!\in\!\F^*_p$ as a root of $\tf_{i-1,\mu}$. So any node $v$  
that has an ancestor at level $\geq\!1$ with $2$ children can have 
no more than $1$ child. Thus, there can be no more than $2(p-1)$ 
nodes at depth $i\!\geq\!2$. It is then clear that $\cT_{p,k}(f)$ has at most 
$1+\left(2\floor{\frac{k-1}{2}} -1\right)(p-1)$ nodes. 

Now, note that the 
coefficient of $x^i$ in the monomial term expansion of $c(\mu+px)^a$ mod 
$p^j$ is simply $c\binom{a}{i}\mu^{a-i}p^i$ mod $p^j$. Since $f$ is a 
trinomial, and Lemma \ref{lem:nodal} 
tell us that $f_{i,\zeta}(x)\!=\!p^{-s}f(\mu+px)$ mod $p^j$ for suitable 
$(s,\mu,j)$, we can then clearly compute the coefficients of 
$x^0,\ldots,x^4$ of any non-root nodal polynomial mod $p^k$  
using $O(\log p)$ multiplications and $O(1)$ additions. 
This takes time $O(k\log^2(p)\log(k\log p))$ via fast modular arithmetic 
\cite{vzgbook}, provided we use Harvey and van der Hoeven's recent fast 
multiplication algorithm \cite{harveymult}. Summing over all 
non-root nodal polynomials, and noting that the cost of reduction mod $p$ 
is negligible compared to the complexity of our earlier steps, we are done. 
\qed 

We can now outline the algorithm that proves Theorem \ref{thm:big}.\\ 
\mbox{}\scalebox{.95}[1]{\fbox{\mbox{}\hspace{.3cm}\vbox{
\begin{algor} {\em
\label{algor:trinosolqp}
{\bf (Solving Trinomial Equations Over $\pmb{\Q^*_p}$)}
\mbox{}\\
{\bf Input.} A prime $p$ and
$c_1,c_2,c_3,a_2,a_3\!\in\!\Z\setminus\{0\}$ with $|c_i|\!\leq\!H$ for all 
$i$ and $1\!\leq\!a_2\!<\!a_3\!=:\!d$. \\
{\bf Output.} \scalebox{.95}[1]{A true declaration that 
$f(x)\!:=\!c_1+c_2x^{a_2}+c_3x^{a_3}$ 
has no roots in $\Q_p$, or $z_1,\ldots,z_m\!\in\!\Q$}\linebreak 
\mbox{}\hspace{1.8cm}\scalebox{.92}[1]{with logarithmic
height $O\!\left(p^2 \log^8(dH)\right)$ such that $m$ is the number of 
roots of $f$ in $\Q_p$, $z_j$}\linebreak  
\mbox{}\hspace{1.8cm}\scalebox{.92}[1]{is an approximate root of 
$f$ with associated true root $\zeta_j\!\in\!\Q_p$ for all $j$, 
and $\#\{\zeta_j\}\!=\!m$.}\\
{\bf Description.} \\
1: If $\ord_p c_1\!\neq\!\ord_p c_2$ mod $a_2$ and $\ord_p c_2\!\neq\!
\ord_p c_3$ mod $a_3-a_2$ then say  
{\tt ``No roots in}\\ 
\mbox{}\hspace{.5cm}{\tt $\Q_p$!''} and {\tt STOP}. \\
2: Rescale and invert roots if necessary, so that we may assume  
$p\nmid c_1 c_2$ and $\ord_p c_3\!\geq\!0$. \\  
3: Compute the mod $p$ reductions of all the nodal polynomials of  
$\cT_{p,k}(f)$, for $k\!:=\!2D+1$\\ 
\mbox{}\hspace{.5cm}where $D\geq 
\max\limits_{\zeta\in\Z_p \; : \; \ord_p \zeta=0}
\ord_p f'(\zeta)$.\\  
4: Use Hensel Lifting to find the first $2D+1$ base-$p$ digits of 
all the non-degenerate roots of\\ 
\mbox{}\hspace{.5cm}$f$ in $\Z_p$ of valuation $0$.\\ 
5: Via Algorithm \ref{algor:binoqp}, or its $p\!=\!2$ version (Algorithm 
\ref{algor:binoq2} from Section \ref{sub:binoq2} of the Appendix),\\  
\mbox{}\hspace{.5cm}find the first 
$O(\log(dH))$ base-$p$ digits of all the degenerate roots of $f$.\\ 
6: If $p|c_3$ then rescale and invert roots to compute approximants for  
the remaining roots of $f$\\ 
\mbox{}\hspace{.5cm}\scalebox{.95}[1]{in $\Q_p$, by computing roots of 
valuation $0$ for a rescaled version of $f$ with coefficients reversed.}} 
\end{algor}}
}} 

\noindent 
{\bf Proof of Theorem \ref{thm:big}:} 
First note that $0$ can not be a root since  
$f(0)\!\neq\!0$ by assumption. So we can focus on roots in $\Q^*_p$. 

\scalebox{.95}[1]{The height bound for our approximate roots from Assertion (1) 
follows directly from Step 3.}  

Assertion (2) follows easily from Theorem \ref{thm:air}: Steps 
3 and 4 (which use Hensel's Lemma) imply a decay rate of $O(p^{-(2D+2^i)})$ 
for the $p$-adic distance of the $i$th Newton iterate to a true root. So the 
decay rate is no worse than $p^{-O(2^i/(2D+1))}$, and thus 
Assertion (2) holds with $\mu\!=\!p^{1/O(p^2\log^8(dH))}$. 

\scalebox{.96}[1]{Assertion (3) on correctly counting the roots of $f$ in 
$\Q_p$ follows immediately from Steps 3--5.}  

So all that remains is to prove correctness (including elaborating 
Step 5) and to do a sufficiently good complexity analysis. 

\smallskip 
\vbox{
\noindent 
{\bf Correctness:} Thanks to Theorem \ref{thm:newt}, Step 1 merely 
guarantees that $f$ has roots of integral valuation, which is a 
necessary condition for their to be roots in $\Q_p$. 
Step 2 merely involves simple substitutions that only negligibly affect 
the heights of the coefficients, similar to the binomial case. 
Steps 3 and 4 correctly count the number of non-degenerate roots 
of $f$ in $\Z_p$ of valuation $0$, thanks to Lemma \ref{lem:ulift}.}

For Step 5, since $0$ is not a root, we can  
rearrange the equations $f(\zeta)\!=\!\zeta f'(\zeta)\!=\!0$ to 
obtain 

\noindent 
that $\zeta\!\in\!\Q^*_p$ is a degenerate root of $f$ if and only 
if $[c_1,c_2\zeta^{a_2},c_3\zeta^{a_3}]^T$ is a right null-vector for
$B\!:=\!\text{\scalebox{1}[.7]{$\begin{bmatrix}1 & 1 & 1\\ 
0 & a_2 & a_3 \end{bmatrix}$}}$. Since $[a_3-a_2,-a_3,a_2]^T$ generates 
the right null-space of $B$ we must have 
$(a_3-a_2)c_2\zeta^{a_2}\!=\!-c_1a_3$ and 
$-a_3c_3\zeta^{a_3-a_2}\!=\!c_2a_2$. 
Via an application of the Extended Euclidean Algorithm, we can then 
find $R,S\!\in\!\Z$ with $Ra_2+S(a_3-a_2)\!=\!\gcd(a_2,a_3)$ and 
the logarithmic heights of $R$ and $S$ of order $O(\log d)$. So 
by multiplying and dividing suitable powers of our binomial equations, 
we get that $\zeta$ must satisfy the single equation $((a_3-a_2)c_2)^R
(-a_3c_3)^S \zeta^{\gcd(a_2,a_3)}\!=\!(-c_1a_3)^R(c_2a_2)^S$. 
The latter equation can then be solved easily, within our overall time 
bound, via Algorithm \ref{algor:binoqp}. Note in particular that 
while the coefficient heights look much larger, any root 
$\zeta$ ultimately satisifies the {\em original} pair of binomials, 
thus implying $\zeta$ must have low logarithmic height. 

\scalebox{.95}[1]{Step 6 merely takes care of the remaining roots, at 
negligible affect to the coefficient heights.}   

Note that we do need to renormalize the roots at the end, due to 
the various rescalings, but this adds a neglible summand of $O(\log H)$ 
to the logarithmic heights of the roots. So we are done.  \qed 

\smallskip 
\noindent 
{\bf Complexity Analysis:} 
Thanks to Lemma \ref{lem:low} (and Step 2 of Algorithm \ref{algor:trinosolqp}) 
our underlying tree $\cT_{p,k}(f)$ will have all its nodal polynomials 
satisfying $\deg\tf_{i,\zeta}\!\leq\!4$ for all $i\!\geq\!1$.   
Steps 3--4 then dominate the overall complexity:  
Theorem \ref{thm:air} tells us that we can take $D\!=\!O\!\left(
p^2\log^8(dH)\right)$, and thus Lemma \ref{lem:complexity} (combined with the 
known upper bounds on the number of $p$-adic rational roots of a trinomial 
\cite{len99,ak11}) implies that the complexity 
of Steps 3--4 is no worse than $O\!\left(p^5\log^3(p)\log^{16}(dH)
\log(p\log(dH)) \right)$, 
assuming we employ brute-force search to find the roots in $\F_p$ of the mod 
$p$ reductions of the nodal polynomials. \qed

\section{Acknowledgements}   
We thank Erich Bach and Bjorn Poonen for informative discussions
on Hensel's Lemma. We also thank the anonymous referees for helpful 
suggestions that improved our paper. 

\bibliographystyle{plain}
\bibliography{20210606arxiv}

\section*{Appendix} 
\subsection{Proof of Proposition \ref{prop:bi}} 
\label{sub:propbi} 
The case $p\!=\!\infty$ follows from an estimate for the 
distance between the vertices of a regular $d$-gon. In particular, 
the minimal spacing between distinct complex roots can easily be expressed 
explicitly as $|c_1/c_2|^{1/d}\sqrt{2(1-\cos\frac{2\pi}{d})}$, which is 
clearly bounded from below by $H^{-1/d} \sqrt{2(1-\cos\frac{2\pi}{d})}$. 
From the elementary 
inequality $1-\cos x\!\geq\!x^2\left(\frac{1}{2!}-\frac{\pi^2}{48}\right)$ 
we easily get $\left|\frac{1}{2}\log\left(1-\cos\frac{2\pi}{d}\right)
\right|\!\leq\!\log(d)-\frac{1}{2}
\log\left(4\pi^2-\frac{\pi^2}{6}\right)$ for all $d\!\geq\!6$. Observing that 
$|\frac{1}{2}\log(1-\cos\frac{2\pi}{d})|\!\leq\!\log 2$ for 
$d\!\in\!\{2,\ldots,5\}$ we get our stated bound via the Triangle Inequality 
applied to 
$\left|\log\left(H^{-1/d}\sqrt{2(1-\cos\frac{2\pi}{d})} \right)\right|$.   

The case of prime $p$ follows easily from the Ultrametric Inequality and 
classical facts on the spacing of $p$-adic roots of unity (see, e.g., 
\cite[Cor.\ 1, Pg.\ 105, Sec.\ 4.3 \& Thm.\ Pg.\ 107, 
Sec.\ 4.4]{Rob00}). 
In particular, when $\gcd(d,p-1)$, the $d$th roots of 
unity in $\C_p$ are all at unit distance. 
At the opposite extreme of 
$d\!=\!p^j$ for $j\!\geq\!1$, the set of distances between distinct $d$th roots 
is exactly $\left\{p^{\frac{-1}{p-1}}, p^{\frac{-1}{p^1(p-1)}}, \ldots,
p^{\frac{-1}{p^{j-1}(p-1)}}\right\}$. So the minimum distance is $p^{-1/(p-1)}$ 
for $d$ a non-trivial $p$th power. 

In complete generality, we see that there are distinct $d$th roots of unity at 
distance $1$ if and only if $d$ is divisible by a prime other than $p$. 
Observing that $\ord_p\!\left(H^{-1/d}\right)\!=\!-\frac{1}{d}\ord_p H\!
\geq\!-\frac{\log H}{d\log p}$ and $|x|_p\!=\!p^{-\ord_p x}$, we then see that 
$\log|H^{-1/d}|_p\!\geq\!-\frac{1}{d}\log H$ and our bound follows again from 
the Triangle Inequality. \qed 

\subsection{Proof of Corollary \ref{cor:binomod}} 
\label{sub:oddbinomod} 
First note that since we are following the notation and assumptions of 
Lemma \ref{lem:binoqp}, we assume $p$ is odd. (We extend Lemma 
\ref{lem:binoqp} and Corollary \ref{cor:binomod} to $p\!=\!2$ in the next 
two sections.) Then $(\Z/p^{2\ell+1})^*$ is cyclic and 
Lemma \ref{lem:binoqp} tells us that we can reduce deciding the feasibility of 
$c_1+c_2x^d\!=\!0$ over $\Q^*_p$ to checking $d\stackrel{?}{|}\ord_p(c_1/c_2)$ 
and $(-c_1/c_2)^r\!\stackrel{?}{=}\!1$ mod $p^{2\ell+1}$ with $r\!=\!
p^\ell(p-1)/\gcd(d,p-1)$.   

The $p$-adic valuation can be computed easily by bisection, ultimately 
resulting in $O(\log H)$ divisions involving integers with 
$O(\max\{\log p,\log H\})\!=\!O(\log(pH))$ bits, and then 
checking divisibility by $d$ involves division by an integer with 
$O(\log d)$ bits. These initial steps dominate 
the computation of the mod $p^{2\ell+1}$ reduction of $-c_1/c_2$. 
From \cite[pp.\ 102--103]{bs} we then see that the $r$th power can be computed 
via recursive squaring using 
just $O(\log(p^{2\ell+1}))\!=\!O(\ell \log p)$ multiplications. 
Since $\ell\!=\!\ord_p d\!\leq\!\log_p d$ we get $\ell\log p\!\leq\!\log d$ and 
thus computing the $r$th power can be done within  
$O(\log(dH))$ multiplications in $\Z/(p^{2\ell+1})$. Each such 
multiplication takes time $O(\ell \log(p) \log (\ell \log p))\!=\!O(\log(d)
\log\log d)$ by fast modular arithmetic 
\cite{vzgbook}, assuming we use the fast integer multiplication algorithm 
of Harvey and van der Hoeven \cite{harveymult}. 
A simple over-estimate of the total complexity then yields our stated 
complexity bound. 

\scalebox{.94}[1]{The remainder of the lemma then follows easily from 
Hensel's Lemma and Proposition \ref{prop:bi}. \qed}  

\subsection{The $2$-adic Version of Lemma \ref{lem:binoqp}}
\label{sub:2adic} 
Recall that the only roots of unity in $\Q_2$ are $\{\pm 1\}$ (see, e.g.,
\cite{Rob00}). The following lemma is then a simple consequence of the
multiplicative group $(\Z/(2^k))^*$ being exactly the product
$\{\pm 1\}\times \left\{1,5,\ldots,5^{2^{k-3}} \text{ mod } 2^k 
\right\}$ (having cardinality $2^{k-1}$) when $k\!\geq\!3$
(see, e.g., \cite[Thm.\ 5.6.2, pg.\ 109 \& Ex.\ 38, pg.\ 192]{bs}),
and Hensel's Lemma.
\begin{lem} {\em \label{lem:binoq2}
Suppose $f(x)\!:=\!c_1+c_2x^d\!\in\!\Z[x]$ with
$|c_1|,|c_2|\!\leq\!H$, and $c_1c_2\!\neq\!0$. Then the number of roots of
the binomial $f$ in $\Q_2$ is either $0$ or $\gcd(d,2)$. In particular,
if $\ell\!:=\!\ord_2 d$ and $u\!:=\!\ord_2(c_2/c_1)$, then $f$ has roots in
$\Q_2$ if and only if {\em both} of the following conditions hold:
(1) $d|u$ and (2) either (i) $d$ is odd or (ii) both
$\frac{c_1}{c_2}2^u\!=\!-1$ mod $8$ and $\left(-\frac{c_1}{c_2}2^u 
\right)^{2^{\ell-1}}\!=\!1$ mod $2^{2\ell+1}$. \qed }
\end{lem}

\subsection{Extending Corollary \ref{cor:binomod} to $\pmb{p\!=\!2}$}  
\label{sub:2binomod}
All three assertions of Corollary \ref{cor:binomod} remain true if we 
replace the assumption that $p$ be odd with the assumption $p\!=\!2$. 
The proof is almost identical to our proof from Section \ref{sub:oddbinomod} 
above, save that use Lemma \ref{lem:binoq2} in place of Lemma \ref{lem:binoqp}. 
In particular, the case $\ell\!=\!0$ remains unchanged. 

As for the case $\ell\!\geq\!1$, the only change is an extra congruence 
condition (mod $8$) to check. However, this additional complexity is 
negligible compared to the other steps, so we are done. \qed 

\subsection{Proof of the $\pmb{p\!=\!2}$ Case of Theorem \ref{thm:binoqp}}  
\label{sub:binoq2} 
Let us first outline an algorithm that proves the $p\!=\!2$ case of Theorem 
\ref{thm:binoqp}.\\ 
\mbox{}\scalebox{.96}[.91]{\fbox{\mbox{}\hspace{.1cm}\vbox{
\begin{algor} {\em
\label{algor:binoq2}
{\bf (Solving Binomial Equations Over $\pmb{\Q^*_2}$)}
\mbox{}\\
{\bf Input.}\hspace{.2cm}$c_1,c_2,d\!\in\!\Z\setminus\{0\}$ with 
$|c_i|\!\leq\!H$ for all $i$. \\ 
{\bf Output.} A true declaration that $f(x)\!:=\!c_1+c_2x^d$ has
no roots in $\Q_2$, or $z_1,$
$\ldots,z_\gamma\!\in\!\Q$ with\\
\mbox{}\hspace{1.8cm}\scalebox{.98}[1]{logarithmic
height $O\!\left(\log\left(dH^{1/d}\right)\right)$ such that $\gamma\!=\!
\gcd(d,2)$, $z_j$ is an approximate}\\
\mbox{}\hspace{1.8cm}\scalebox{.96}[1]{root of $f$ with associated true root
$\zeta_j\!\in\!\Q_p$ for all $j$, and the $\zeta_j$ are pair-wise distinct.}\\
{\bf Description.} \\
1: \scalebox{.95}[1]{If $\ord_2 c_1\!\neq\!\ord_2 c_2$ mod $d$ then say
{\tt ``No roots in $\Q_p$!''} and {\tt STOP}.} \\
2: Let $\ell\!:=\!\ord_2 d$ and replace $f$ with $f(x)\!:=\!c'_1+c'_2x^d$ where 
$c'_i\!:=\!\frac{c_i}{2^{\ord_2 c_i}}$ for all $i$. \\
3: \scalebox{.93}[1]{If $c'_1\!\neq\!-c'_2$ mod $8$ or $\left(-\frac{c'_1}{c'_2}
\right)^{2^{\ell-1}}\!\!\!\!\!\!\!\!\!\neq\!1$
mod $2^{2\ell+1}$ then say {\tt ``No roots in $\Q_2$!''} and {\tt STOP}.}\\
4: \scalebox{.95}[1]{Let $\delta\!:=\!1$. If $d\!\leq\!-1$ then
set $\delta\!:=\!-1$ and respectively replace $d$ by $|d|$ and $f(x)$ by 
$x^df(1/x)$.}\\
5: Let $x_1\!:=\!1$. If $\gamma\!=\!1$ then {\tt GOTO} Step 7.\\  
6: Let $x_2\!:=\!3$. \\ 
7: Output $\left\{x_1 2^{\ord_2(c_1/c_2)/d},\ldots,x_\gamma 
2^{\ord_2(c_1/c_2)/d}\right\}$. }
\end{algor}}
}}

Similar to the case of odd $p$, it clearly suffices to prove the correctness of 
Algorithm \ref{algor:binoq2}, and then analyze its complexity. 

\noindent 
{\bf Correctness:} The proof is almost the same 
as the Correctness proof for odd $p$, save that we respectively 
replace Lemma \ref{lem:binoqp} and Algorithm \ref{algor:binoqp} by 
Lemma \ref{lem:binoq2} and Algorithm \ref{algor:binoq2}. 
In particular, Steps 5--8 of Algorithm \ref{algor:binoqp} collapse into 
Steps 5--6 of Algorithm \ref{algor:binoq2}. 

So we must explain Steps 5--6 here: These steps merely give us the mod 
$4$ reductions of the $\gamma$ many roots of $f$ in $\Z_2$, since 
Steps 5 and 6 are executed only after Steps 1 and 3 certify that 
$f$ indeed has roots in $\Z_2$. (Remember 
that $\gamma\!\in\!\{1,2\}$ for $p\!=\!2$.) Furthermore, Hensel's Lemma 
implies that the root $1$ of $\tf$ lifts to the sole root of 
$f$ in $\Z_2$ when $\ell\!=\!0$. So the case $\ell\!=\!0$ is done. 

If $\ell\!\geq\!1$ then there is one more complication: The nodal 
polynomial $\tf_{1,1}$ is now quadratic. This is because 
Lemma \ref{lem:binodepth} tells us that $\deg \tf_{1,1}\!\leq\!2$. 
Furthermore, $\ell\!\geq\!1$ implies that $\gamma\!=\!2$ and thus $f$ must 
have exactly $2$ roots in $\Z_2$. So then, Lemma \ref{lem:ulift} tells 
us that $\deg \tf_{1,1}\!\leq\!1$ would imply that $f$ has $\leq\!1$ 
root in $\Z_2$. Therefore, $\tf_{1,1}$ is quadratic. 

Furthermore, $\tf_{1,1}$ must also have $2$ distinct roots: This is 
because $\tf_{1,1}$ equal to $x^2$ or $1+x^2\!=\!(1+x)^2$ mod $2$ would 
imply that no nodal polynomial $\tf_{i,\zeta}$, for $i\!\geq\!1$, has a 
non-degenerate root. So, again by Lemma \ref{lem:ulift}, we would 
not attain $2$ roots in $\Z_2$. (Similarly, it is impossible for 
$\tf_{1,1}$ to be irreducible.) Therefore, the mod $4$ reductions 
of the two roots of $f$ in $\Z_2$ must be $1$ and $3$. So Steps 5--6 
are indeed correct.  

Lemma \ref{lem:nodal} then tells us that Hensel's Lemma
--- applied to $f_{1,1}(x)\!=\!2^{-s(f,1)}f(1+2x)$ and {\em either} 
start point $0$ or $1$ in $\Z/(2)$ --- implies that $1+0$ and $1+1\cdot 2$ are  
approximate roots of $f$ with distinct associated true roots in $\Z_2$. So
Steps 5--7 indeed give us suitable approximants in $\Q$ to all the roots
of $f$ in $\Q_2$, and our algorithm is correct.

Note also that the outputs, being integers in $\{1,3\}$ rescaled by a 
factor of $2^{\ord_2(c_1/c_2)/d}$ (or possibly the reciprocals
of such quantities), clearly each have bit-length\\ 
\mbox{}\hfill $O\!\left(\frac{|\log(c_1/c_2)|}{d\log 2}\log 2\right)\!=\!
O\!\left(\frac{\log H}{d}\right)\!=\!O\!\left(\log\!\left(
H^{1/d}\right)\right)$. \hfill \qed

\medskip 
\noindent 
{\bf Complexity Analysis:} We merely use the same techniques 
as for Algorithm \ref{algor:binoqp}, save for Steps 5--8 there 
being collapsed into Steps 5--6 here. Also, the prime $p$ is fixed 
to $2$. So we easily arrive at an overall complexity bound of 
$O([\log(dH)\log\log(dH)]^2)$. \qed 

\subsection{The Proof of Case 2 of Theorem \ref{thm:tri} (Both Roots Large)} 
\label{sub:bigroots} 
Simply observe that $1/\zeta_1$ and $1/\zeta_2$ are roots of        
the reciprocal polynomial $g(x)\!:=\!x^{\deg f}f(\frac{1}{x})$.   
In particular, we can apply Case 1 to the trinomial $g$ since  
$\norm{\frac{1}{\zeta_1}}, \norm{\frac{1}{\zeta_2}}<1$. 
We then obtain $\norm{\frac{1}{\zeta_1} - \frac{1}{\zeta_2}} 
\geq e^{-O(M)}$. Hence $\norm{\zeta_1-\zeta_2} = \norm{\zeta_1}
\norm{\zeta_2} \norm{\frac{1}{\zeta_1} - \frac{1}{\zeta_2}}
\geq \norm{\frac{1}{\zeta_1} - \frac{1}{\zeta_2}}\geq 
e^{-O(M)}$. \qed 

\subsection{The Proof of the $p\!=\!2$ Case of Theorem \ref{thm:tetra} } 
Returning to where we observed that $G$ is square-free when 
$p$ is odd, assume instead that $p\!=\!2$. Then, as $h\!>\!2$, we 
have $p^{d(h-1)/2+1} \geq 8$.
Then, as $G(x) = 1-x^2  = (3-x)(5-x) \mod 2^3$, we obtain that  
$G$ is square-free in $\Z_2[x]$. 
Hensel's Lemma then implies that there are 
roots $\zeta_1,\zeta_2 \in \Z_p$ of $G$ such that 
$\zeta_1 \equiv 3 \mod p^{d(h-1)/2+1}$ and $\zeta_2 \equiv 5 \mod 
p^{d(d-1)/2+1}$. We then proceed as in the remainder of the proof 
of the case of odd $p$. \qed 

\subsection{The Proof of the $p\!=\!\infty$ Case of Theorem 
\ref{thm:tetra}} \label{sub:infinity}
Shifting by $\frac{1}{2^{h-1}}$, we get\\  
$g(x):=f_{d,\frac{1}{2}}(x+2^{1-h}) 
= (x+2^{1-h})^d - 2^{2h}x^2$\\ 
\mbox{}\hspace{4.1cm}\scalebox{.85}[1]{$=2^{d(1-h)} + d2^{(d-1)(1-h)}x 
+\left({d\choose 2} 2^{(d-2)(1-h)} -2^{2h}\right)x^2 + {d\choose 3} 
2^{(d-3)(1-h)} x^3 + \cdots+ x^d$.}\linebreak   
We will see momentarily that, unlike $\anewt(f)$ (which has $3$ lower
edges), $\anewt(g)$ will have just $2$ lower edges.
(See the right-hand illustration in Example \ref{ex:newts}.)
This will force (via Theorem \ref{thm:newt}) the existence of two distinct
roots of small norm for $g$, thus yielding two nearby roots of $f$ after
undoing our earlier shift.

Toward this end, note that the three lowest order terms of $g$ contribute
the points\linebreak 
\scalebox{.91}[1]{$p_0:= (0,d(h-1)\log2)$, $p_1:=(1,(d-1)(h-1)\log2-\log d)$,
and $p_2 = \left(2,-\log\left(4^h - \frac{{d\choose 2}}
{2^{(d-2)(h-1)}} \right)\right)$}\linebreak 
as potential vertices of $\anewt(g)$.
Observe that $\frac{{d\choose 2}}{2^{(d-2)(h-1)}}\!<\!0.059$ for all 
$h\!\geq\!3$ and $d\!\geq\!4$, and thus $p_2$ is the only point of $\anewt(f)$ 
with negative $y$-coordinate. So $p_2$ is a vertex of $\anewt(f)$, and
all edges with vertices to the right of $p_2$ have positive slope.
Furthermore, the slopes of the line segments $\overline{p_0p_1}$
and $\overline{p_0p_2}$ are respectively $-(h-1)\log(2)-\log d$ and
a number less than $-\frac{1}{2}\log(4^h-0.059)-\frac{1}{2}d(h-1)\log 2$.

Since $2^{h-1}\!<\!\sqrt{4^h-0.059}$ and $\log d\!<\!\frac{1}{2}d(h-1)\log 2$
for all $d\!\geq\!4$ and $h\!\geq\!3$, we thus see that the slope
of $\overline{p_0p_2}$ is more negative. So the leftmost lower edge of
$\anewt(g)$ has vertices $p_0$ and $p_2$. It is easily checked that
the slope of this edge is less than $-10.3$, which is in turn clearly
$<\!-2\log 3$. So by Theorem \ref{thm:newt}, there are
two roots $z_1,z_2$ of $g$ such that
\begin{align*}
\log|z_i| \leq \frac{1}{2}\left[-\log\left(2^{2h} - {d\choose 2}
2^{(d-2)(1-h)}\right)  - d(h-1)\log2\right].
\end{align*}
These two roots thus satisfy $|z_i| = 2^{-\Omega(dh)}$. Now, for
$i\!\in\!\{1,2\}$, $\zeta_i = z_i+2^{1-h}$ yields roots of $f_{d,\frac{1}{2}}$
with\\ 
$|\zeta_1-\zeta_2| = |z_1+2^{1-h}-(z_2+2^{1-h})| \leq |z_1|+|z_2| < 
2^{-\Omega(dh)}$. \qed

\subsection{The Proof of Theorem \ref{thm:air}} \label{sub:gcd} 
Suppose $r\!:=\!\gcd(a_2,a_3)$. The special case $r\!=\!1$ was 
proved in an even stronger form (summing over roots in $\C_p$ instead of 
roots in $\Z_p$) in \cite[Sec.\ 5]{airr}. So assume $r\!>\!1$. We can then 
clearly write $f(x)\!=\!g(x^r)$ for some trinomial $g\!\in\Z[x]$ with 
relatively prime exponents and coefficients satisfying the conditions of 
Theorem \ref{thm:air}. 

By the Chain Rule we have $f'(x)\!=\!x^{r-1}g'(x^r)$. 
So then our desired sum of valuations equals the sum of 
$\ord_p g'(\zeta^r)$ over the roots $\zeta\!\in\!\Z_p$ of 
$f$ with $|\zeta|_p\!=\!1$. Since the $r$th power map induces a 
$\gcd(r,p-1)$-to-$1$ endomorphism of $\Q^*_p$, we see that 
our desired sum is then at most 
$\gcd(r,p-1)O(p\log^8(dH))\!=\!O(p^2\log^8(dH))$. \qed 
\mbox{}\\

\end{document}